\documentclass[12pt,a4paper, twoside,reqno]{amsart}

\usepackage{amscd}
\usepackage{amsmath,amssymb,amsthm,mathrsfs}
\usepackage[margin=2cm]{geometry}

\usepackage{pgf,tikz} % for drawing vector graphics inside latex
\usetikzlibrary{matrix,arrows}

\usepackage{hyperref} % just for better pdf navigation

\usepackage{parskip}
\usepackage{enumerate}
\DeclareMathOperator{\shiftoperator}{\Delta}
\newcommand{\shift}[1]{\shiftoperator_{#1}}

\DeclareMathOperator{\cro}{cr} %cross ratio
\DeclareMathOperator{\sr}{sr} %star ratio
\DeclareMathOperator{\mr}{mr} %multi ratio

\newcommand{\obs}{\mathcal{G}} %surface graph class

\DeclareMathOperator{\miq}{miq} %miquel dynamics
\DeclareMathOperator{\cli}{cli} %clifford dynamics
\DeclareMathOperator{\mut}{mut} %mutation of combinatorics
\DeclareMathOperator{\mob}{mob} %möbius mutation map

\newcommand{\Z}{\mathbb{Z}} % whole numbers
\newcommand{\R}{\mathbb{R}} % real numbers
\newcommand{\C}{\mathbb{C}} % complex numbers

\newcommand{\Ce}{\hat{\mathbb{C}}} % the riemann sphere
\newcommand{\CP}{\hat{\mathbb{C}}_{\circ}} % the set of circle patterns

\newcommand{\f}[2]{\frac{#1}{#2}} % fractions fractions

\newcommand{\menge}[2][]{ % a set in brackets
	  \if\relax\detokenize{#1}\relax
		  \left\{#2\right\}
	  \else
		  \left\{#1\ | \ #2\right\}
	  \fi}

\newcommand{\ggborientationexample}{
	\definecolor{ududff}{rgb}{0.30196078431372547,0.30196078431372547,1.}
	\definecolor{xdxdff}{rgb}{0.49019607843137253,0.49019607843137253,1.}
	\definecolor{yqyqyq}{rgb}{0.5019607843137255,0.5019607843137255,0.5019607843137255}
	\definecolor{ffffff}{rgb}{1.,1.,1.}
	\definecolor{cqcqcq}{rgb}{0.7529411764705882,0.7529411764705882,0.7529411764705882}
	\begin{tikzpicture}[line cap=round,line join=round,>=triangle 45,x=1.0cm,y=1.0cm]
	\clip(-0.6571102672898803,-6.627659181405768) rectangle (11.181168334608852,1.756032782216452);
	\draw [line width=2.pt,color=cqcqcq] (3.3942754777070054,-3.0606827229299354) circle (1.7783954347257358cm);
	\draw [line width=2.pt,color=cqcqcq] (6.77546181095711,-2.0243413769402183) circle (2.1370942848404817cm);
	\draw [line width=2.pt,color=cqcqcq] (3.7426762411900922,0.12886083169559615) circle (2.132889190322997cm);
	\draw [line width=2.pt,color=cqcqcq] (0.9821131574203864,-1.4361726129015606) circle (1.4631361069847129cm);
	\draw [line width=2.pt,color=cqcqcq] (0.14772288535075676,-4.133365747680886) circle (2.0448748996090926cm);
	\draw [line width=2.pt,color=cqcqcq] (6.012095462000406,-4.9793078890306885) circle (1.8040964229209795cm);
	\draw [-latex,line width=2.pt] (1.6449281854232318,-2.7405665358466074) -- (2.44,-1.56);
	\draw [-latex,line width=2.pt] (4.650027305998744,-1.801406078186893) -- (2.44,-1.56);
	\draw [-latex,line width=2.pt] (1.6245460771436373,-0.12162051761653891) -- (2.44,-1.56);
	\draw [-latex,line width=2.pt] (4.650027305998744,-1.801406078186893) -- (5.14,-3.4);
	\draw [-latex,line width=2.pt] (4.650027305998744,-1.801406078186893) -- (5.864174395728069,-0.09127975631271701);
	\draw [-latex,line width=2.pt] (1.6449281854232318,-2.7405665358466074) -- (2.18,-4.36);
	\draw [-latex,line width=2.pt] (4.243489823263154,-4.623220824082986) -- (2.18,-4.36);
	\draw [-latex,line width=2.pt] (4.243489823263154,-4.623220824082986) -- (5.14,-3.4);
	\draw [-latex,line width=2.pt] (7.54,-4.02) -- (5.14,-3.4);
	\draw [-latex,line width=2.pt] (1.6449281854232318,-2.7405665358466074) -- (-0.3014719811126467,-2.138437847481427);
	\draw [line width=2.pt,color=cqcqcq] (-0.6465885198379875,0.11920317421649092) circle (2.2838669856528373cm);
	\draw [-latex,line width=2.pt] (1.6245460771436373,-0.12162051761653891) -- (-0.30147198111264695,-2.138437847481427);
	\draw [-latex,line width=2.pt] (1.6245460771436373,-0.12162051761653891) -- (1.6234643278375969,0.37001879756349776);
	\draw [line width=2.pt,color=cqcqcq] (6.732071119297279,0.6944783687112227) circle (1.1707521316757477cm);
	\draw [line width=2.pt,color=cqcqcq] (3.09827029101989,-5.381181841834758) circle (1.3733290506877007cm);
	\draw [line width=2.pt,color=cqcqcq] (0.9677886698958211,-5.401866129418487) circle (1.5984183871760238cm);
	\draw [-latex,line width=2.pt] (-0.4798568675078866,-6.0795561311873865) -- (2.18,-4.36);
	\draw [-latex,line width=2.pt] (2.2,-6.42) -- (2.18,-4.36);
	\draw [line width=2.pt,color=cqcqcq] (8.132237301796751,-5.1283038211574805) circle (1.2566154469971191cm);
	\draw [line width=2.pt,color=cqcqcq] (8.85813953488372,-4.69813953488372) circle (1.4823511939468932cm);
	\draw [line width=2.pt,color=cqcqcq] (9.174110139050008,-3.3157681416562084) circle (1.7793983412525696cm);
	\draw [line width=2.pt,color=cqcqcq] (10.027520652094806,-0.87481476880145) circle (1.3536684285130531cm);
	\draw [line width=2.pt,color=cqcqcq] (9.18103691698583,-0.04025453058853846) circle (1.5566059209530023cm);
	\draw [-latex,line width=2.pt] (5.637100901716218,1.108846292703647) -- (5.864174395728069,-0.09127975631271701);
	\draw [-latex,line width=2.pt] (7.624599892224223,-0.06318444766610853) -- (5.864174395728069,-0.09127975631271701);
	\draw [-latex,line width=2.pt] (7.624599892224223,-0.06318444766610853) -- (8.862363279500968,-1.5638913330838415);
	\draw [-latex,line width=2.pt] (7.624599892224223,-0.06318444766610853) -- (7.894276487190005,0.8356847306948038);
	\draw [-latex,line width=2.pt] (10.7,0.3) -- (8.862363279500968,-1.5638913330838415);
	\draw [-latex,line width=2.pt] (10.7,0.3) -- (7.894276487190005,0.8356847306948039);
	\draw [-latex,line width=2.pt] (10.509506438192993,-2.139768489833295) -- (8.862363279500968,-1.5638913330838415);
	\draw [-latex,line width=2.pt] (7.54,-4.02) -- (8.862363279500968,-1.5638913330838413);
	\draw [-latex,line width=2.pt] (7.54,-4.02) -- (10.34,-4.66);
	\draw [-latex,line width=2.pt] (7.54,-4.02) -- (8.82,-6.18);
	\draw [-latex,line width=2.pt] (7.54,-4.02) -- (7.390811457117948,-6.142882599998282);
	\draw [-latex,line width=2.pt] (4.243489823263154,-4.623220824082986) -- (4.405907495389701,-5.800845544352244);
	\draw [-latex,line width=2.pt] (2.2,-6.42) -- (4.405907495389701,-5.800845544352244);
	\draw [-latex,line width=2.pt] (-0.4798568675078866,-6.0795561311873865) -- (-0.30147198111264684,-2.138437847481427);
	\draw [-latex,line width=2.pt] (5.637100901716218,1.108846292703647) -- (1.6234643278375973,0.37001879756349776);
	\draw [-latex,line width=2.pt] (5.637100901716218,1.108846292703647) -- (7.894276487190005,0.8356847306948038);
	\draw [-latex,line width=2.pt] (5.793888290362814,-6.770159511427752) -- (4.405907495389701,-5.800845544352244);
	\draw [-latex,line width=2.pt] (5.793888290362814,-6.770159511427752) -- (7.390811457117948,-6.142882599998282);
	\draw [-latex,line width=2.pt,color=cqcqcq] (3.7426762411900922,0.12886083169559615) -- (6.77546181095711,-2.0243413769402183);
	\draw [-latex,line width=2.pt,color=cqcqcq] (6.77546181095711,-2.0243413769402183) -- (6.732071119297279,0.6944783687112226);
	\draw [-latex,line width=2.pt,color=cqcqcq] (6.732071119297279,0.6944783687112227) -- (3.7426762411900922,0.12886083169559615);
	\draw [-latex,line width=2.pt,color=cqcqcq] (6.77546181095711,-2.0243413769402183) -- (3.3942754777070054,-3.0606827229299354);
	\draw [-latex,line width=2.pt,color=cqcqcq] (3.3942754777070054,-3.0606827229299354) -- (3.7426762411900922,0.12886083169559637);
	\draw [-latex,line width=2.pt,color=cqcqcq] (3.3942754777070054,-3.0606827229299354) -- (6.012095462000406,-4.9793078890306885);
	\draw [-latex,line width=2.pt,color=cqcqcq] (6.012095462000406,-4.9793078890306885) -- (6.77546181095711,-2.0243413769402183);
	\draw [-latex,line width=2.pt,color=cqcqcq] (6.77546181095711,-2.0243413769402183) -- (9.174110139050008,-3.3157681416562084);
	\draw [-latex,line width=2.pt,color=cqcqcq] (6.012095462000406,-4.9793078890306885) -- (3.09827029101989,-5.381181841834758);
	\draw [-latex,line width=2.pt,color=cqcqcq] (3.09827029101989,-5.381181841834758) -- (3.3942754777070054,-3.0606827229299354);
	\draw [-latex,line width=2.pt,color=cqcqcq] (3.3942754777070054,-3.0606827229299354) -- (0.14772288535075662,-4.133365747680886);
	\draw [-latex,line width=2.pt,color=cqcqcq] (0.14772288535075676,-4.133365747680886) -- (0.9821131574203864,-1.4361726129015606);
	\draw [-latex,line width=2.pt,color=cqcqcq] (0.9821131574203864,-1.4361726129015606) -- (3.3942754777070054,-3.0606827229299354);
	\draw [-latex,line width=2.pt,color=cqcqcq] (0.9821131574203864,-1.4361726129015606) -- (-0.6465885198379874,0.11920317421649096);
	\draw [-latex,line width=2.pt,color=cqcqcq] (3.7426762411900922,0.12886083169559615) -- (0.9821131574203865,-1.4361726129015606);
	\draw [-latex,line width=2.pt,color=cqcqcq] (-0.6465885198379875,0.11920317421649092) -- (3.7426762411900922,0.12886083169559615);
	\draw [-latex,line width=2.pt,color=cqcqcq] (6.732071119297279,0.6944783687112227) -- (9.18103691698583,-0.04025453058853845);
	\draw [-latex,line width=2.pt,color=cqcqcq] (9.18103691698583,-0.04025453058853846) -- (6.77546181095711,-2.0243413769402183);
	\draw [line width=2.pt,color=cqcqcq] (3.8531698127965934,-8.070718351123332) circle (2.336202368091699cm);
	\draw [line width=2.pt,color=cqcqcq] (2.0592331392374703,-8.362730073644283) circle (1.9478232589562272cm);
	\draw [line width=2.pt,color=cqcqcq] (-0.40421161325195054,-7.193830148815722) circle (1.1168387488142757cm);
	\draw [line width=2.pt,color=cqcqcq] (-3.052607770269616,-3.988510765075995) circle (3.3153458237108864cm);
	\draw [line width=2.pt,color=cqcqcq] (7.009487089832101,-7.518470048183641) circle (1.427462469100697cm);
	\draw [line width=2.pt,color=cqcqcq] (8.491110922145793,-7.2706875744755095) circle (1.1391959491883132cm);
	\draw [line width=2.pt,color=cqcqcq] (9.65563365646902,-6.7496447497432825) circle (1.0113252437934446cm);
	\draw [line width=2.pt,color=cqcqcq] (11.316365218034196,-5.705467831487235) circle (1.4304866394558136cm);
	\draw [line width=2.pt,color=cqcqcq] (11.046872250154777,-3.441726901300143) circle (1.4084948423954664cm);
	\draw [line width=2.pt,color=cqcqcq] (11.395853675443252,-2.0952245979350144) circle (0.8874658209117763cm);
	\draw [line width=2.pt,color=cqcqcq] (11.881027297142653,-0.09005128365061922) circle (1.2437706703703777cm);
	\draw [line width=2.pt,color=cqcqcq] (9.554096106129196,1.9136950197220826) circle (1.9791683482372395cm);
	\draw [line width=2.pt,color=cqcqcq] (6.890220350952134,2.001289553738856) circle (1.538428850555164cm);
	\draw [line width=2.pt,color=cqcqcq] (3.5596222633420114,1.1232907159619572) circle (2.0775288528115357cm);
	\draw [-latex,line width=2.pt,color=cqcqcq] (-3.052607770269616,-3.988510765075995) -- (0.14772288535075662,-4.133365747680886);
	\draw [-latex,line width=2.pt,color=cqcqcq] (0.14772288535075676,-4.133365747680886) -- (0.9677886698958211,-5.401866129418487);
	\draw [-latex,line width=2.pt,color=cqcqcq] (0.9677886698958211,-5.401866129418487) -- (3.09827029101989,-5.381181841834758);
	\draw [-latex,line width=2.pt,color=cqcqcq] (3.09827029101989,-5.381181841834758) -- (3.8531698127965934,-8.070718351123332);
	\draw [-latex,line width=2.pt,color=cqcqcq] (3.8531698127965934,-8.070718351123332) -- (2.0592331392374703,-8.362730073644283);
	\draw [-latex,line width=2.pt,color=cqcqcq] (2.0592331392374703,-8.362730073644283) -- (0.9677886698958211,-5.401866129418487);
	\draw [-latex,line width=2.pt,color=cqcqcq] (0.9677886698958211,-5.401866129418487) -- (-0.40421161325195065,-7.193830148815722);
	\draw [-latex,line width=2.pt,color=cqcqcq] (3.8531698127965934,-8.070718351123332) -- (6.012095462000406,-4.9793078890306885);
	\draw [-latex,line width=2.pt,color=cqcqcq] (6.012095462000406,-4.9793078890306885) -- (7.009487089832101,-7.518470048183641);
	\draw [-latex,line width=2.pt,color=cqcqcq] (7.009487089832101,-7.518470048183641) -- (8.132237301796751,-5.1283038211574805);
	\draw [-latex,line width=2.pt,color=cqcqcq] (8.132237301796751,-5.1283038211574805) -- (8.491110922145793,-7.2706875744755095);
	\draw [-latex,line width=2.pt,color=cqcqcq] (11.046872250154777,-3.441726901300143) -- (9.174110139050008,-3.3157681416562084);
	\draw [-latex,line width=2.pt,color=cqcqcq] (9.18103691698583,-0.04025453058853846) -- (9.554096106129196,1.9136950197220826);
	\draw [-latex,line width=2.pt,color=cqcqcq] (6.890220350952134,2.001289553738856) -- (6.732071119297279,0.6944783687112226);
	\draw [-latex,line width=2.pt,color=cqcqcq] (3.7426762411900922,0.12886083169559615) -- (3.5596222633420114,1.1232907159619572);
	\draw [-latex,line width=2.pt] (-2.7845506176341415,-0.6840193830751009) -- (-0.30147198111264695,-2.138437847481427);
	\draw [-latex,line width=2.pt] (-0.4798568675078866,-6.0795561311873865) -- (-1.484229881881109,-6.909419652736003);
	\draw [-latex,line width=2.pt] (-0.4798568675078866,-6.0795561311873865) -- (0.6912128961823197,-6.976174606768837);
	\draw [-latex,line width=2.pt] (2.2,-6.42) -- (0.6912128961823198,-6.976174606768837);
	\draw [-latex,line width=2.pt] (0.12,-8.18) -- (0.6912128961823198,-6.976174606768837);
	\draw [-latex,line width=2.pt] (5.793888290362814,-6.770159511427752) -- (6.120183742147127,-8.635066956869743);
	\draw [-latex,line width=2.pt] (7.824743670867615,-6.346716743187516) -- (7.390811457117948,-6.142882599998282);
	\draw [-latex,line width=2.pt] (10.509506438192993,-2.139768489833295) -- (10.34,-4.66);
	\draw [-latex,line width=2.pt] (10.509506438192993,-2.139768489833295) -- (11.33914386548157,-1.2095725326516322);
	\draw [-latex,line width=2.pt] (10.7,0.3) -- (11.33914386548157,-1.2095725326516322);
	\draw [-latex,line width=2.pt] (10.7,0.3) -- (11.32,1.02);
	\draw [-latex,line width=2.pt] (5.637100901716218,1.108846292703647) -- (5.36,2.16);
	\draw [-latex,line width=2.pt] (12.02,-4.46) -- (10.34,-4.66);
	\draw [-latex,line width=2.pt] (9.887124973750588,-5.76517000717387) -- (10.34,-4.66);
	\draw [-latex,line width=2.pt] (9.887124973750588,-5.76517000717387) -- (8.82,-6.18);
	\draw [-latex,line width=2.pt] (9.887124973750588,-5.76517000717387) -- (10.643181347905578,-6.96765462537324);
	\draw [-latex,line width=2.pt] (7.824743670867615,-6.346716743187516) -- (8.82,-6.18);
	\draw [-latex,line width=2.pt] (9.523494873444708,-7.7523002692287735) -- (8.82,-6.18);
	\draw [-latex,line width=2.pt] (7.824743670867615,-6.346716743187516) -- (8.161639927175043,-8.361199506606487);
	\draw [-latex,line width=2.pt] (2.2,-6.42) -- (2.8088712487379945,-10.160522651478372);
	\draw [-latex,line width=2.pt] (7.968680452407601,3.0984163184242974) -- (7.894276487190005,0.8356847306948039);
	\draw [-latex,line width=2.pt] (11.851173098323466,-1.3334636075747592) -- (11.33914386548157,-1.2095725326516322);
	\draw [-latex,line width=2.pt] (10.509506438192993,-2.139768489833295) -- (12.148983947668995,-2.564682132056885);
	\draw [-latex,line width=2.pt] (-0.1832165500371627,2.3555696079115327) -- (1.6234643278375969,0.37001879756349787);
	\draw [-latex,line width=2.pt] (2.4731096846141014,2.894057936838389) -- (1.6234643278375969,0.37001879756349787);
	\draw [line width=2.pt,color=cqcqcq] (1.2947823447009326,1.8856840920283398) circle (1.5508943003582953cm);
	\draw [-latex,line width=2.pt,color=cqcqcq] (1.2947823447009326,1.8856840920283398) -- (-0.6465885198379875,0.11920317421649096);
	\draw [-latex,line width=2.pt,color=cqcqcq] (3.5596222633420114,1.1232907159619572) -- (1.2947823447009323,1.8856840920283398);
	\draw [-latex,line width=2.pt,color=cqcqcq] (3.5596222633420114,1.1232907159619572) -- (6.890220350952134,2.001289553738856);
	\draw [-latex,line width=2.pt,color=cqcqcq] (9.554096106129196,1.9136950197220826) -- (6.890220350952134,2.001289553738856);
	\draw [-latex,line width=2.pt,color=cqcqcq] (9.554096106129196,1.9136950197220826) -- (11.881027297142653,-0.0900512836506191);
	\draw [-latex,line width=2.pt,color=cqcqcq] (9.174110139050008,-3.3157681416562084) -- (10.027520652094806,-0.8748147688014503);
	\draw [-latex,line width=2.pt,color=cqcqcq] (10.027520652094806,-0.87481476880145) -- (9.18103691698583,-0.04025453058853845);
	\draw [-latex,line width=2.pt,color=cqcqcq] (11.881027297142653,-0.09005128365061922) -- (10.027520652094806,-0.87481476880145);
	\draw [-latex,line width=2.pt,color=cqcqcq] (11.395853675443252,-2.0952245979350144) -- (11.881027297142653,-0.09005128365061932);
	\draw [-latex,line width=2.pt,color=cqcqcq] (10.027520652094806,-0.87481476880145) -- (11.395853675443252,-2.0952245979350144);
	\draw [-latex,line width=2.pt,color=cqcqcq] (11.395853675443252,-2.0952245979350144) -- (11.046872250154777,-3.441726901300143);
	\draw [-latex,line width=2.pt,color=cqcqcq] (9.174110139050008,-3.3157681416562084) -- (8.85813953488372,-4.69813953488372);
	\draw [-latex,line width=2.pt,color=cqcqcq] (8.85813953488372,-4.69813953488372) -- (8.132237301796751,-5.1283038211574805);
	\draw [-latex,line width=2.pt,color=cqcqcq] (8.132237301796751,-5.1283038211574805) -- (6.012095462000406,-4.9793078890306885);
	\draw [-latex,line width=2.pt,color=cqcqcq] (7.009487089832101,-7.518470048183641) -- (3.8531698127965934,-8.070718351123332);
	\draw [-latex,line width=2.pt,color=cqcqcq] (8.491110922145793,-7.2706875744755095) -- (7.009487089832101,-7.518470048183641);
	\draw [-latex,line width=2.pt,color=cqcqcq] (8.491110922145793,-7.2706875744755095) -- (9.65563365646902,-6.7496447497432825);
	\draw [-latex,line width=2.pt,color=cqcqcq] (9.65563365646902,-6.7496447497432825) -- (8.85813953488372,-4.69813953488372);
	\draw [-latex,line width=2.pt,color=cqcqcq] (11.316365218034196,-5.705467831487235) -- (9.65563365646902,-6.7496447497432825);
	\draw [-latex,line width=2.pt,color=cqcqcq] (8.85813953488372,-4.69813953488372) -- (11.316365218034196,-5.705467831487235);
	\draw [-latex,line width=2.pt,color=cqcqcq] (11.316365218034196,-5.705467831487235) -- (11.046872250154777,-3.441726901300143);
	\draw [-latex,line width=2.pt,color=cqcqcq] (-0.40421161325195054,-7.193830148815722) -- (2.0592331392374703,-8.362730073644283);
	\draw [-latex,line width=2.pt,color=cqcqcq] (-0.40421161325195054,-7.193830148815722) -- (-3.052607770269616,-3.988510765075995);
	\draw [-latex,line width=2.pt,color=cqcqcq] (-0.6465885198379875,0.11920317421649092) -- (-3.052607770269616,-3.9885107650759952);
	\begin{scriptsize}
	\draw [fill=black] (2.44,-1.56) circle (2.5pt);
	\draw [fill=black] (2.18,-4.36) circle (2.5pt);
	\draw [fill=black] (5.14,-3.4) circle (2.5pt);
	\draw [fill=ffffff] (4.650027305998744,-1.801406078186893) circle (2.5pt);
	\draw [fill=ffffff] (7.54,-4.02) circle (2.5pt);
	\draw [fill=black] (5.864174395728069,-0.09127975631271701) circle (2.5pt);
	\draw [fill=ffffff] (1.6449281854232318,-2.7405665358466074) circle (2.5pt);
	\draw [fill=ffffff] (1.6245460771436373,-0.12162051761653891) circle (2.5pt);
	\draw [fill=black] (-0.30147198111264684,-2.138437847481427) circle (2.5pt);
	\draw [fill=ffffff] (4.243489823263154,-4.623220824082986) circle (2.5pt);
	\draw [fill=black] (1.6234643278375969,0.37001879756349776) circle (2.5pt);
	\draw [fill=black] (8.862363279500968,-1.5638913330838415) circle (2.5pt);
	\draw [fill=ffffff] (7.624599892224223,-0.06318444766610853) circle (2.5pt);
	\draw [fill=ffffff] (5.637100901716218,1.108846292703647) circle (2.5pt);
	\draw [fill=ffffff] (2.2,-6.42) circle (2.5pt);
	\draw [fill=ffffff] (-0.4798568675078866,-6.0795561311873865) circle (2.5pt);
	\draw [fill=black] (7.390811457117948,-6.142882599998282) circle (2.5pt);
	\draw [fill=black] (8.82,-6.18) circle (2.5pt);
	\draw [fill=black] (10.34,-4.66) circle (2.5pt);
	\draw [fill=ffffff] (10.509506438192993,-2.139768489833295) circle (2.5pt);
	\draw [fill=ffffff] (10.7,0.3) circle (2.5pt);
	\draw [fill=black] (7.894276487190005,0.8356847306948038) circle (2.5pt);
	\draw [fill=black] (4.405907495389701,-5.800845544352244) circle (2.5pt);
	\draw [fill=ffffff] (5.793888290362814,-6.770159511427752) circle (2.5pt);
	\draw [fill=yqyqyq] (-0.6465885198379875,0.11920317421649092) circle (2.5pt);
	\draw [fill=yqyqyq] (0.9821131574203864,-1.4361726129015606) circle (2.5pt);
	\draw [fill=yqyqyq] (3.7426762411900922,0.12886083169559615) circle (2.5pt);
	\draw [fill=yqyqyq] (6.732071119297279,0.6944783687112227) circle (2.5pt);
	\draw [fill=yqyqyq] (9.18103691698583,-0.04025453058853846) circle (2.5pt);
	\draw [fill=yqyqyq] (10.027520652094806,-0.87481476880145) circle (2.5pt);
	\draw [fill=yqyqyq] (6.77546181095711,-2.0243413769402183) circle (2.5pt);
	\draw [fill=yqyqyq] (9.174110139050008,-3.3157681416562084) circle (2.5pt);
	\draw [fill=yqyqyq] (8.85813953488372,-4.69813953488372) circle (2.5pt);
	\draw [fill=yqyqyq] (8.132237301796751,-5.1283038211574805) circle (2.5pt);
	\draw [fill=yqyqyq] (6.012095462000406,-4.9793078890306885) circle (2.5pt);
	\draw [fill=yqyqyq] (3.09827029101989,-5.381181841834758) circle (2.5pt);
	\draw [fill=yqyqyq] (3.3942754777070054,-3.0606827229299354) circle (2.5pt);
	\draw [fill=yqyqyq] (0.14772288535075676,-4.133365747680886) circle (2.5pt);
	\draw [fill=xdxdff] (0.6912128961823198,-6.976174606768837) circle (2.5pt);
	\draw [fill=ududff] (0.12,-8.18) circle (2.5pt);
	\draw [fill=xdxdff] (-1.484229881881109,-6.909419652736003) circle (2.5pt);
	\draw [fill=xdxdff] (6.120183742147127,-8.635066956869743) circle (2.5pt);
	\draw [fill=ffffff] (7.824743670867615,-6.346716743187516) circle (2.5pt);
	\draw [fill=xdxdff] (8.161639927175043,-8.361199506606487) circle (2.5pt);
	\draw [fill=ffffff] (9.887124973750588,-5.76517000717387) circle (2.5pt);
	\draw [fill=xdxdff] (9.523494873444708,-7.7523002692287735) circle (2.5pt);
	\draw [fill=ududff] (12.02,-4.46) circle (2.5pt);
	\draw [fill=xdxdff] (11.33914386548157,-1.2095725326516322) circle (2.5pt);
	\draw [fill=xdxdff] (12.148983947668995,-2.564682132056885) circle (2.5pt);
	\draw [fill=ududff] (11.32,1.02) circle (2.5pt);
	\draw [fill=ududff] (5.36,2.16) circle (2.5pt);
	\draw [fill=yqyqyq] (-0.40421161325195054,-7.193830148815722) circle (2.5pt);
	\draw [fill=yqyqyq] (-3.052607770269616,-3.988510765075995) circle (2.5pt);
	\draw [fill=yqyqyq] (3.5596222633420114,1.1232907159619572) circle (2.5pt);
	\draw [fill=yqyqyq] (6.890220350952134,2.001289553738856) circle (2.5pt);
	\draw [fill=yqyqyq] (9.554096106129196,1.9136950197220826) circle (2.5pt);
	\draw [fill=yqyqyq] (11.881027297142653,-0.09005128365061922) circle (2.5pt);
	\draw [fill=yqyqyq] (11.395853675443252,-2.0952245979350144) circle (2.5pt);
	\draw [fill=yqyqyq] (11.046872250154777,-3.441726901300143) circle (2.5pt);
	\draw [fill=yqyqyq] (11.316365218034196,-5.705467831487235) circle (2.5pt);
	\draw [fill=yqyqyq] (9.65563365646902,-6.7496447497432825) circle (2.5pt);
	\draw [fill=yqyqyq] (8.491110922145793,-7.2706875744755095) circle (2.5pt);
	\draw [fill=yqyqyq] (7.009487089832101,-7.518470048183641) circle (2.5pt);
	\draw [fill=yqyqyq] (3.8531698127965934,-8.070718351123332) circle (2.5pt);
	\draw [fill=yqyqyq] (2.0592331392374703,-8.362730073644283) circle (2.5pt);
	\draw [fill=yqyqyq] (0.9677886698958211,-5.401866129418487) circle (2.5pt);
	\draw [fill=yqyqyq] (-2.7845506176341415,-0.6840193830751009) circle (2.5pt);
	\draw [fill=yqyqyq] (10.643181347905578,-6.96765462537324) circle (2.5pt);
	\draw [fill=yqyqyq] (2.8088712487379945,-10.160522651478372) circle (2.5pt);
	\draw [fill=yqyqyq] (7.968680452407601,3.0984163184242974) circle (2.5pt);
	\draw [fill=yqyqyq] (11.851173098323466,-1.3334636075747592) circle (2.5pt);
	\draw [fill=xdxdff] (-0.1832165500371627,2.3555696079115327) circle (2.5pt);
	\draw [fill=xdxdff] (2.4731096846141014,2.894057936838389) circle (2.5pt);
	\draw [fill=yqyqyq] (1.2947823447009326,1.8856840920283398) circle (2.5pt);
	\end{scriptsize}
	\end{tikzpicture}
}

\newcommand{\ggbcliffordcomb}[5]{ %  
	\definecolor{yqyqyq}{rgb}{0.75,0.75,0.75}
	\begin{tikzpicture}[line cap=round,line join=round,>=triangle 45,x=1.0cm,y=1.0cm,scale=.5,line width=1.5pt,scale=1.1]
	\clip(-1.7,-0.38) rectangle (4.86,5.2);
	\draw [color=yqyqyq] (-0.012684563758389415,3.4862416107382543) circle (1.5935056506165621cm);
	\draw [color=yqyqyq] (2.7982587533364742,1.6902103940964046) circle (1.9604285635211869cm);
	\draw [color=yqyqyq] (0.5791070720167777,1.6023166417049615) circle (1.8786732543194797cm);
	\draw [color=yqyqyq] (2.961360901520605,4.2929121098596115) circle (1.7759395077687445cm);

	\ifnum#5=0
	\tikzstyle{dira}=[-latex];
	\else
	\tikzstyle{dira}=[latex-];
	\fi
		
	\ifnum#5=0		
		\draw [dira, color=black] (-1.14,2.36) -- (0.9,2.18);
		\draw [dira] (4.22,3.04) -- (0.9,2.18);
		\draw [dira] (0.9,2.18) -- (1.2,4.52);
		\draw [dira] (0.9,2.18) -- (1.68,0.08);
		\draw [fill=black] (0.9,2.18) circle (2.5pt);
	\else
		\draw [dira, color=black]  (2.05,2.77) -- (-1.14,2.36) ;
		\draw [dira]  (2.05,2.77) -- (4.22,3.04);
		\draw [dira] (1.2,4.52) --  (2.05,2.77);
		\draw [dira] (1.68,0.08) --  (2.05,2.77);	
		\draw [fill=black] (2.05,2.77) circle (2.5pt);
	\fi

	\draw [fill=black] (1.2,4.52) circle (2.5pt);
	\draw [fill=black] (-1.14,2.36) circle (2.5pt);
	\draw [fill=black] (1.68,0.08) circle (2.5pt);
	\draw [fill=black] (4.22,3.04) circle (2.5pt);

	\ifnum#1=1
	\draw [dira] (1.2,4.52) -- (4.22,3.04);
	\fi
	\ifnum#2=1
	\draw [dira] (4.22,3.04) -- (1.68,0.08);
	\fi
	\ifnum#3=1
	\draw [dira] (1.68,0.08) -- (-1.14,2.36);
	\fi
	\ifnum#4=1
	\draw [dira] (-1.14,2.36) -- (1.2,4.52);
	\fi

	\end{tikzpicture}
}

\newcommand{\ggbmiqcomb}[5]{ %  which legs and then which central circle
	\definecolor{yqyqyq}{rgb}{0.75,0.75,0.75}
	\begin{tikzpicture}[line cap=round,line join=round,>=triangle 45,x=1.0cm,y=1.0cm,scale=.7, line width=1.5pt]
	\clip(-2.22,1.1) rectangle (2.02,5.66);
	\draw [color=yqyqyq] (-0.13170583377130393,2.9623582851565433) circle (1.6935932842182857cm);
	\draw [color=yqyqyq] (2.4061709806041622,3.138049366815499) circle (1.8022926687675231cm);
	\draw [color=yqyqyq] (-0.10585348912684864,5.806116196046681) circle (1.9057276141591373cm);
	\draw [color=yqyqyq] (-0.1250931981095718,1.2944126595544845) circle (1.386741828734134cm);
	\draw [color=yqyqyq] (-2.809028340080972,3.1614979757085018) circle (1.9317654523655317cm);
	\draw [color=yqyqyq] (-0.16970092170609585,3.8098767458921814) circle (1.7630576280226038cm);
	\ifnum#5=0	
		\draw [line width=2.pt] (-1.22,4.26)-- (0.98,4.24);
		\draw [line width=2.pt] (0.98,4.24)-- (1.1454475908339077,1.8500915265197229);
		\draw [line width=2.pt] (1.1454475908339077,1.8500915265197229)-- (-1.4,1.84);
		\draw [line width=2.pt] (-1.4,1.84)-- (-1.22,4.26);	
	\else
		\draw [line width=2.pt] (-1.6759951585796644,4.726091716131247)-- (-1.08,2.3);
		\draw [line width=2.pt] (-1.08,2.3)-- (0.7935594599130316,2.333224573939806);
		\draw [line width=2.pt] (0.7935594599130316,2.333224573939806)-- (1.39206296608269,4.627963390770062);
		\draw [line width=2.pt] (1.39206296608269,4.627963390770062)-- (-1.6759951585796644,4.726091716131247);	
	\fi
	\ifnum#1=1
		\draw [line width=2.pt] (-1.22,4.26)-- (-1.6759951585796644,4.726091716131247);
	\fi
	\ifnum#2=1
		\draw [line width=2.pt] (0.98,4.24)-- (1.39206296608269,4.627963390770062);
	\fi
	\ifnum#3=1
		\draw [line width=2.pt] (0.7935594599130316,2.333224573939806)-- (1.1454475908339077,1.8500915265197229);
	\fi
	\ifnum#4=1
		\draw [line width=2.pt] (-1.08,2.3)-- (-1.4,1.84);
	\fi
	\end{tikzpicture}
}

\newcommand{\ggbmutationcomb}[5]{ %  which legs and then which orientation
	\begin{tikzpicture}[line cap=round,line join=round,>=triangle 45,x=1.0cm,y=1.0cm,scale=.8,line width=1.5pt]
	\clip(-2,-2) rectangle (2,2);
	\definecolor{gray}{rgb}{0.75,0.75,0.75}
	
	\ifnum#5=0
		\tikzstyle{dira}=[-latex];
		\definecolor{col0}{rgb}{0,0,0}
		\definecolor{col1}{rgb}{1,1,1}	
	\else
		\tikzstyle{dira}=[latex-];
		\definecolor{col1}{rgb}{0,0,0}
		\definecolor{col0}{rgb}{1,1,1}	
	\fi

	\coordinate (f) at (0,0);
	\coordinate (v1) at (1,0);
	\coordinate (v2) at (0,1);
	\coordinate (v3) at (-1,0);
	\coordinate (v4) at (0,-1);
	\coordinate (f12) at (1.5,1.5);
	\coordinate (f23) at (-1.5,1.5);
	\coordinate (f34) at (-1.5,-1.5);
	\coordinate (f14) at (1.5,-1.5);
	\coordinate (c1) at (1.85,0);
	\coordinate (c2) at (0,1.85);
	\coordinate (c3) at (-1.85,0);
	\coordinate (c4) at (0,-1.85);

	\draw [-, black] (v1)-- (v2) -- (v3) -- (v4) -- (v1);
	\draw [dira,gray] (f)-- (f12);
	\draw [dira,gray] (f23)-- (f);
	\draw [dira,gray] (f)-- (f34);
	\draw [dira,gray] (f14)-- (f);
	
	\def\rad{3.5pt}	
	\ifnum#1=1
		\draw [-, black] (v1) -- (c1);
		\draw [dira,gray] (f12) -- (f14);
		\draw [fill=col0] (c1) circle (\rad);
	\fi
	\ifnum#2=1
		\draw [-, black] (v2) -- (c2);
		\draw [dira,gray] (f12) -- (f23);
		\draw [fill=col1] (c2) circle (\rad);
	\fi
	\ifnum#3=1
		\draw [-, black] (v3) -- (c3);
		\draw [dira,gray] (f34) -- (f23);
		\draw [fill=col0] (c3) circle (\rad);
	\fi
	\ifnum#4=1
		\draw [-,black] (v4) -- (c4);
		\draw [dira,gray] (f34) -- (f14);
		\draw [fill=col1] (c4) circle (\rad);
	\fi

	\draw [fill=col1] (v1) circle (\rad);
	\draw [fill=col0] (v2) circle (\rad);
	\draw [fill=col1] (v3) circle (\rad);	
	\draw [fill=col0] (v4) circle (\rad);	
		
	\end{tikzpicture}
}

\newcommand{\ggbmiquelmenelaus}{
	\definecolor{yqyqyq}{rgb}{0.75,0.75,0.75}
	\definecolor{ffffff}{rgb}{1.,1.,1.}
	\begin{tikzpicture}[line cap=round,line join=round,>=triangle 45,x=1.0cm,y=1.0cm,scale=.8, line width=1.2pt]
	\clip(-2.26,-1.76) rectangle (6.3,5.5);
	\draw [color=yqyqyq] (2.36,0.42) circle (2.08153789300123cm);
	\draw [domain=-2.26:6.3] plot(\x,{(--5.82215106988162-2.1443502733353323*\x)/1.8130581543100863});
	\draw [domain=-2.26:6.3] plot(\x,{(-3.8868683403670596--1.9064786372407228*\x)/1.4581458179072524});
	\draw [color=yqyqyq] (0.6738287988556445,2.4142778269238385) circle (1.7682502658584827cm);
	\draw [color=yqyqyq] (-0.3757264693898448,3.655613699682519) circle (3.0436423202006706cm);
	\draw [color=yqyqyq] (4.6561069295090896,3.4220857695920133) circle (2.4002858349584746cm);
	\draw [color=yqyqyq] (1.383811966290962,-0.8563343757124051) circle (3.5185954021653214cm);
	\draw [color=yqyqyq] (0.41229922467052843,3.6190413423321495) circle (2.3159930829171924cm);
	\draw [domain=-2.26:6.3] plot(\x,{(--2.8669763202044076-0.036572357350369344*\x)/0.7880256940603731});
	\draw [domain=-2.26:6.3] plot(\x,{(--1.4432094045326447-1.204763515408311*\x)/0.2615295741851161});
	\begin{normalsize}
	\draw [fill=ffffff] (2.36,0.42) circle (3.4pt);
	\draw[color=black] (2.9,0.5) node {$M$};
	\draw [fill=black] (2.44,2.5) circle (3.4pt);
	\draw[color=black] (3.,2.8) node {$V_1$};
	\draw [fill=black] (4.346478637240723,1.0418541820927476) circle (3.4pt);
	\draw[color=black] (4.6,1.46) node {$V_4$};
	\draw [fill=black] (0.2956497266646676,0.6869418456899137) circle (3.4pt);
	\draw[color=black] (-.2,0.2) node {$V_2$};
	\draw [fill=ffffff] (0.6738287988556445,2.4142778269238385) circle (3.4pt);
	\draw[color=black] (0.2,2.0) node {$M_{23}$};
	\draw [fill=ffffff] (-0.3757264693898448,3.655613699682519) circle (3.4pt);
	\draw[color=black] (-1.2,4) node {$M_{12}$};
	\draw [fill=ffffff] (4.6561069295090896,3.4220857695920133) circle (3.4pt);
	\draw[color=black] (4.25,3.9) node {$M_{14}$};
	\draw [fill=ffffff] (1.383811966290962,-0.8563343757124051) circle (3.4pt);
	\draw[color=black] (2,-1) node {$M_{34}$};
	\draw [fill=black] (-0.9689191018153591,1.7599932271130188) circle (3.4pt);
	\draw[color=black] (-1.8,1.9) node {$V_{234}$};
	\draw [fill=black] (2.5349302666785065,4.545465337932655) circle (3.4pt);
	\draw[color=black] (1.8,4.7) node {$V_{124}$};
	\draw [fill=ffffff] (0.41229922467052843,3.6190413423321495) circle (3.4pt);
	\draw[color=black] (1.3,3.2) node {$M_{1234}$};
	\end{normalsize}
	\end{tikzpicture}	
}

\newcommand{\ggbmiquel}{	
	\definecolor{ffffff}{rgb}{1.,1.,1.}
	\begin{tikzpicture}[line cap=round,line join=round,>=triangle 45,x=1.0cm,y=1.0cm,scale=.52, line width=1.2pt]
	\clip(-4.86,-5.18) rectangle (6.76,6.3);
	\draw [color=lightgray] (0.98,0.08) circle (2.3820159529272678cm);
	\draw [color=lightgray] (3.02,3.18) circle (2.6953292934259445cm);
	\draw [color=lightgray] (4.16,-1.78) circle (2.4301061856756045cm);
	\draw [color=lightgray] (-2.04,3.46) circle (2.697035409482049cm);
	\draw [color=lightgray] (-1.2089316455933286,-2.0940365648580483) circle (2.971332342867912cm);
	\draw [color=lightgray] (0.9567261437753695,1.2475875050907475) circle (3.020941156711098cm);
	\begin{normalsize}
	\draw [fill=ffffff] (0.98,0.08) circle (5pt);
	\draw[color=black] (1.65,0.) node {$M$};
	\draw [fill=black] (0.44,2.4) circle (5pt);
	\draw[color=black] (-0.2,2.9) node {$V_{1}$};
	\draw [fill=ffffff] (3.02,3.18) circle (5pt);
	\draw[color=black] (4,3.6) node {$M_{14}$};
	\draw [fill=ffffff] (4.16,-1.78) circle (5pt);
	\draw[color=black] (4.8,-2.5) node {$M_{34}$};
	\draw [fill=black] (3.324353306805309,0.5019094368119901) circle (5pt);
	\draw[color=black] (2.6,1.) node {$V_{4}$};
	\draw [fill=ffffff] (-2.04,3.46) circle (5pt);
	\draw[color=black] (-1.8,4) node {$M_{12}$};
	\draw [fill=black] (-1.2648066664070712,0.8767703749853977) circle (5pt);
	\draw[color=black] (-.8,0.2) node {$V_{2}$};
	\draw [fill=black] (1.7614245362366083,-2.170194590289338) circle (5pt);
	\draw[color=black] (2.5,-2.5) node {$V_{3}$};
	\draw [fill=black] (1.732290369331135,-1.6721017941924947) circle (5pt);
	\draw[color=black] (.9,-1.26) node {$V_{234}$};
	\draw [fill=black] (-2.025096455398336,0.7630057685715551) circle (5pt);
	\draw[color=black] (-2.8,1.45) node {$V_{123}$};
	\draw [fill=black] (3.9154778082020596,0.6377727294717268) circle (5pt);
	\draw[color=black] (4.8,0) node {$V_{134}$};
	\draw [fill=black] (0.5418126314149989,4.239899696285337) circle (5pt);
	\draw[color=black] (1.9,4.7) node {$V_{124}$};
	\draw [fill=ffffff] (0.9567261437753695,1.2475875050907475) circle (5pt);
	\draw[color=black] (.6,.9) node {$M_{1234}$};
	\draw [fill=ffffff] (-1.2089316455933286,-2.0940365648580483) circle (5pt);
	\draw[color=black] (-1.7,-2.5) node {$M_{23}$};
	\end{normalsize}
	\end{tikzpicture}
}

\newcommand{\ggbmenelaus}{
	\begin{tikzpicture}[line cap=round,line join=round,>=triangle 45,x=1.0cm,y=1.0cm, scale=1, line width=1.2pt]
	\clip(-2.98,0.) rectangle (5.52,4.8);
	\draw [domain=-2.98:5.52] plot(\x,{(-9.744--1.92*\x)/-2.72});
	\draw [domain=-2.98:5.52] plot(\x,{(-5.6644-3.22*\x)/-1.26});
	\draw [domain=-2.98:5.52] plot(\x,{(--4.2308--1.3*\x)/3.98});
	\draw [domain=-2.98:5.52] plot(\x,{(--13.507400263526344--0.18979287898175246*\x)/5.5315271129705135});
	\begin{large}
	\draw [fill=black] (2.44,1.86) circle (2.5pt);
	\draw [fill=black] (-0.28,3.78) circle (2.5pt);
	\draw[color=black] (2.45,1.66) node[below] {$V_{14}$};
	\draw[color=black] (.6,3.9) node {$V_{1234}$};
	\draw[color=black] (-.9,0.3) node {$V_{12}$};
	\draw[color=black] (-0.55,2.22) node[below] {$V_{23}$};
	\draw[color=black] (1.77,2.92) node {$V_{34}$};
	\draw[color=black] (4.84,2.45) node[below] {$V$};

	\draw[color=black] (-1.87,4.4) node {$c_{134}$};
	\draw [fill=black] (-1.54,0.56) circle (2.5pt);
	\draw[color=black] (0.6,4.6) node {$c_{123}$};
	\draw[color=black] (-2.61,0.6) node {$c_1$};
	\draw [fill=black] (-0.8145430073067556,2.413945647993847) circle (2.5pt);
	\draw [fill=black] (4.716984105663758,2.6037385269755995) circle (2.5pt);
	\draw[color=black] (-2.75,2.7) node {$c_3$};
	\draw [fill=black] (1.5407585407509048,2.4947586771170083) circle (2.5pt);
	\end{large}
	\end{tikzpicture}	
}

\newcommand{\ggbcliffordthree}{
	\definecolor{ffffff}{rgb}{1.,1.,1.}
	\definecolor{yqyqyq}{rgb}{0.75,0.75,0.75}
	\begin{tikzpicture}[line cap=round,line join=round,>=triangle 45,x=1.0cm,y=1.0cm,scale=1.4,line width=1.2pt]
	\clip(-4.5,0.8913726707688842) rectangle (1.4291474129312973,6.105497389025843);
	\draw [color=yqyqyq] (-1.0010603733151935,3.605595248715718) circle (1.489708342285813cm);
	\draw [color=yqyqyq] (-0.2789835298137834,4.34503133136104) circle (1.489096784536228cm);
	\draw [color=yqyqyq] (-2.9357061276435688,4.500580629054133) circle (1.3731756306066816cm);
	\draw [color=yqyqyq] (-1.0883864947168855,2.494171885421456) circle (1.5341725846776013cm);
	\draw [] (-1.7003405069089776,3.9010115374565952)-- (-2.9357061276435688,4.500580629054133);
	\draw [] (-2.9357061276435688,4.500580629054133)-- (-2.44,3.22);
	\draw [] (-2.44,3.22)-- (-1.0883864947168855,2.494171885421456);
	\draw [] (-1.0883864947168855,2.494171885421456)-- (-1.7003405069089776,3.9010115374565952);
	\draw [] (-2.9357061276435688,4.500580629054133)-- (-1.6388124258775196,4.95188766291424);
	\draw [] (-1.6388124258775196,4.95188766291424)-- (-0.2789835298137834,4.34503133136104);
	\draw [] (-0.2789835298137834,4.34503133136104)-- (-1.7003405069089776,3.9010115374565952);
	\draw [] (-0.2789835298137834,4.34503133136104)-- (0.36,3.);
	\draw [] (0.36,3.)-- (-1.0883864947168855,2.494171885421456);
	\draw [] (-1.0010603733151935,3.605595248715718)-- (-1.6388124258775196,4.95188766291424);
	\draw [] (-1.0010603733151935,3.605595248715718)-- (-2.44,3.22);
	\draw [] (-1.0010603733151935,3.605595248715718)-- (0.36,3.);
	\begin{normalsize}
	\draw [fill=black] (-1.6388124258775196,4.95188766291424) circle (2.1pt);
	\draw[color=black] (-2,5.1) node {$V_{12}$};
	\draw [fill=black] (-2.44,3.22) circle (2.1pt);
	\draw[color=black] (-2.8,3.35) node {$V_{23}$};
	\draw [fill=black] (0.36,3.) circle (2.1pt);
	\draw[color=black] (0.7,2.9) node {$V_{13}$};
	\draw [fill=black] (-1.7003405069089776,3.9010115374565952) circle (2.1pt);
	\draw[color=black] (-1.4,3.8) node {$V$};
	\draw [fill=ffffff] (-1.0883864947168855,2.494171885421456) circle (2.1pt);
	\draw[color=black] (-0.9,2.77) node {$M_3$};
	\draw [fill=ffffff] (-2.9357061276435688,4.500580629054133) circle (2.1pt);
	\draw[color=black] (-3.1914386126134056,4.7) node {$M_2$};
	\draw [fill=ffffff] (-1.0010603733151935,3.605595248715718) circle (2.1pt);
	\draw[color=black] (-0.7,3.8) node {$M_{123}$};
	\draw [fill=ffffff] (-0.2789835298137834,4.34503133136104) circle (2.1pt);
	\draw[color=black] (-0.1411005209530001,4.55) node {$M_1$};
	\end{normalsize}
	\end{tikzpicture}	
}
\newcommand{\ggbcliffordfour}{
	\definecolor{ffffff}{rgb}{1.,1.,1.}
	\definecolor{yqyqyq}{rgb}{0.75,0.75,0.75}
	\begin{tikzpicture}[line cap=round,line join=round,>=triangle 45,x=1.0cm,y=1.0cm, line width=1.2pt]
	\clip(-3.38,-1.18) rectangle (4.4,6.14);
	\draw [color=yqyqyq] (-0.3402094647013189,3.31211016291699) circle (1.774543964231919cm);
	\draw [color=yqyqyq] (-1.2200111255330983,1.0318078991285) circle (1.8862402956014401cm);
	\draw [color=yqyqyq] (2.068892928558444,0.8336811488820213) circle (1.80475370011392cm);
	\draw [color=yqyqyq] (1.6515085854456257,3.164077064595257) circle (1.8032592942390715cm);
	\draw [color=yqyqyq] (-0.5080459821279759,1.689836289245356) circle (1.8408826066714092cm);
	\draw [color=yqyqyq] (-0.8428404910489878,4.012747351159195) circle (1.733663845558456cm);
	\draw [color=yqyqyq] (2.3189542597524953,3.8261401124964274) circle (1.8144095604020753cm);
	\draw [] (-0.019777849657817148,2.486917014017243)-- (-0.8428404910489878,4.012747351159195);
	\draw [] (-0.8428404910489878,4.012747351159195)-- (-2.02,2.74);
	\draw [] (-2.02,2.74)-- (-1.2200111255330983,1.0318078991285);
	\draw [] (-1.2200111255330983,1.0318078991285)-- (-0.019777849657817148,2.486917014017243);
	\draw [] (0.42,0.1)-- (-1.2200111255330983,1.0318078991285);
	\draw [] (3.2,2.24)-- (1.6515085854456257,3.164077064595257);
	\draw [] (1.6515085854456257,3.164077064595257)-- (-0.019777849657817148,2.486917014017243);
	\draw [] (1.6515085854456257,3.164077064595257)-- (0.74,4.72);
	\draw [] (0.74,4.72)-- (-0.8428404910489878,4.012747351159195);
	\draw [color=yqyqyq] (1.5268005655096193,1.5378945924688128) circle (1.8145380543993739cm);
	\draw [] (0.42,0.1)-- (1.5268005655096193,1.5378945924688128);
	\draw [] (1.5268005655096193,1.5378945924688128)-- (-0.019777849657817148,2.486917014017243);
	\draw [] (1.5268005655096193,1.5378945924688128)-- (3.2,2.24);
	\draw [] (3.2,2.24)-- (2.068892928558444,0.8336811488820213);
	\draw [] (2.068892928558444,0.8336811488820213)-- (0.42,0.1);
	\draw [] (1.1868780914695725,2.4082241763125127)-- (2.068892928558444,0.8336811488820213);
	\draw [] (1.1868780914695725,2.4082241763125127)-- (2.3189542597524953,3.8261401124964274);
	\draw [] (3.2,2.24)-- (2.3189542597524953,3.8261401124964274);
	\draw [] (1.1868780914695725,2.4082241763125127)-- (-0.3402094647013189,3.31211016291699);
	\draw [] (-0.3402094647013189,3.31211016291699)-- (0.74,4.72);
	\draw [] (0.74,4.72)-- (2.3189542597524953,3.8261401124964274);
	\draw [] (-0.3402094647013189,3.31211016291699)-- (-2.02,2.74);
	\draw [] (-2.02,2.74)-- (-0.5080459821279759,1.689836289245356);
	\draw [] (-0.5080459821279759,1.689836289245356)-- (1.1868780914695725,2.4082241763125127);
	\draw [] (-0.5080459821279759,1.689836289245356)-- (0.42,0.1);
	\draw [] (-0.8428404910489878,4.012747351159195)-- (0.6458178209037198,3.1242166765898554);
	\draw [] (0.6458178209037198,3.1242166765898554)-- (2.3189542597524953,3.8261401124964274);
	\draw [] (0.52,1.76)-- (-1.2200111255330983,1.0318078991285);
	\draw [] (0.52,1.76)-- (2.068892928558444,0.8336811488820213);
	\draw [] (-0.3402094647013189,3.31211016291699)-- (0.52,1.76);
	\draw [] (0.52,1.76)-- (1.6515085854456257,3.164077064595257);
	\draw [] (0.6458178209037198,3.1242166765898554)-- (-0.8428404910489878,4.012747351159195);
	\draw [] (0.6458178209037198,3.1242166765898554)-- (1.5268005655096193,1.5378945924688128);
	\draw [] (0.6458178209037198,3.1242166765898554)-- (-0.5080459821279759,1.689836289245356);
	\begin{normalsize}
	\draw [fill=black] (0.52,1.76) circle (3pt);
	\draw[color=black] (0.75,1.3) node {$V_{13}$};
	\draw [fill=black] (0.74,4.72) circle (3pt);
	\draw[color=black] (0.8,5.14) node {$V_{12}$};
	\draw [fill=black] (-2.02,2.74) circle (3pt);
	\draw[color=black] (-2.6,2.8) node {$V_{23}$};
	\draw [fill=black] (0.42,0.1) circle (3pt);
	\draw[color=black] (0.5,-0.4) node {$V_{34}$};
	\draw [fill=black] (3.2,2.24) circle (3pt);
	\draw[color=black] (3.8,2.2) node {$V_{14}$};
	\draw [fill=black] (1.1868780914695725,2.4082241763125127) circle (3pt);
	\draw[color=black] (1.95,2.5) node {$V_{1234}$};
	\draw [fill=black] (-0.019777849657817148,2.486917014017243) circle (3pt);
	\draw[color=black] (-.5,2.5) node {$V$};
	\draw [fill=ffffff] (-0.8428404910489878,4.012747351159195) circle (3pt);
	\draw[color=black] (-1,4.4) node {$M_2$};
	\draw [fill=ffffff] (-0.3402094647013189,3.31211016291699) circle (3pt);
	\draw[color=black] (-1,3.4) node {$M_{123}$};
	\draw [fill=ffffff] (1.6515085854456257,3.164077064595257) circle (3pt);
	\draw[color=black] (2.3,3.1) node {$M_1$};
	\draw [fill=ffffff] (2.3189542597524953,3.8261401124964274) circle (3pt);
	\draw[color=black] (2.59,4.2) node {$M_{124}$};
	\draw [fill=ffffff] (2.068892928558444,0.8336811488820213) circle (3pt);
	\draw[color=black] (2.6,.7) node {$M_{134}$};
	\draw [fill=ffffff] (-1.2200111255330983,1.0318078991285) circle (3pt);
	\draw[color=black] (-1.5,.7) node {$M_3$};
	\draw [fill=ffffff] (-0.5080459821279759,1.689836289245356) circle (3pt);
	\draw[color=black] (-1.2,1.8) node {$M_{234}$};
	\draw [fill=black] (0.6458178209037198,3.1242166765898554) circle (3pt);
	\draw[color=black] (0.87,3.54) node {$V_{24}$};
	\draw [fill=ffffff] (1.5268005655096193,1.5378945924688128) circle (3pt);
	\draw[color=black] (2.1,1.44) node {$M_4$};
	\end{normalsize}
	\end{tikzpicture}	
}

\newcommand{\ggbsroctahedron}{
	\begin{tikzpicture}[line cap=round,line join=round,>=latex,x=1.0cm,y=1.0cm,line width=1.25pt]
	\clip(-2.26,-0.88) rectangle (2.94,4.14);
	\draw [->] (1.2,1.6) -- (-1.74,3.78);
	\draw [->] (1.2,1.6) -- (2.44,-0.34);
	\draw [->] (2.08,3.84) -- (1.2,1.6);
	\draw [->] (-1.68,-0.52) -- (1.2,1.6);
	\draw [->] (2.44,-0.34) -- (-1.68,-0.52);
	\draw [->] (-1.74,3.78) -- (-1.68,-0.52);
	\draw [->] (-0.6583505712162934,1.8055876791802874) -- (2.44,-0.34);
	\draw [->] (-0.6583505712162934,1.8055876791802874) -- (-1.74,3.78);
	\draw [->] (-1.68,-0.52) -- (-0.6583505712162934,1.8055876791802876);
	\draw [->] (2.08,3.84) -- (-0.6583505712162934,1.8055876791802872);
	\draw [->] (-1.74,3.78) -- (2.08,3.84);
	\draw [->] (2.44,-0.34) -- (2.08,3.84);
	\begin{large}
	\draw [fill=black] (2.08,3.84) circle (3.5pt);
	\draw[color=black] (2.47,3.84) node {$z_{\bar3}$};
	\draw [fill=black] (1.2,1.6) circle (3.5pt);
	\draw[color=black] (1.61,1.7) node {$z_1$};
	\draw [fill=black] (-1.74,3.78) circle (3.5pt);
	\draw[color=black] (-2.1,3.66) node {$z_2$};
	\draw [fill=black] (2.44,-0.34) circle (3.5pt);
	\draw[color=black] (2.8,-.1) node {$z_{\bar2}$};
	\draw [fill=black] (-1.68,-0.52) circle (3.5pt);
	\draw[color=black] (-2.03,-0.2) node {$z_3$};
	\draw [fill=black] (-0.6583505712162934,1.8055876791802874) circle (3.5pt);
	\draw[color=black] (-1.1,1.85) node {$z_{\bar1}$};
	\end{large}
	\end{tikzpicture}
}

\newcommand{\ggblocalintegrability}{
	\definecolor{cqcqcq}{rgb}{0.7529411764705882,0.7529411764705882,0.7529411764705882}
	\definecolor{eqeqeq}{rgb}{0.8784313725490196,0.8784313725490196,0.8784313725490196}
	\begin{tikzpicture}[line cap=round,line join=round,>=triangle 45,x=1.0cm,y=1.0cm,scale=.5]
	\clip(-9.370319320036238,-4.5667934572512445) rectangle (17.738149733093422,8.944683231690323);
	\draw [line width=2.pt,color=cqcqcq] (-0.17738423104116285,0.051053006057898954) circle (2.3187997556127136cm);
	\draw [line width=2.pt] (15.867473943901507,64.04329888210773) circle (63.729929542321976cm);
	\draw [line width=2.pt] (2.339188518206574,-22.364821225267747) circle (24.494392639543495cm);
	\draw [line width=2.pt,color=cqcqcq] (13.461571194742321,0.2414721161567196) circle (2.3187997556127136cm);
	\draw [line width=2.pt,color=cqcqcq] (-6.6238816030720296,2.790818474278457) circle (2.3187997556127136cm);
	\begin{large}
	\draw [fill=eqeqeq] (-0.17738423104116285,0.051053006057898954) circle (5.5pt);
	\draw[] (0.3143309016260143,-0.8) node {$M_{12}$};
	\draw[] (-2.9,-1) node {$q_1$};
	\draw [fill=eqeqeq] (6.578354113429635,4.711657306990274) circle (5.5pt);
	\draw[] (7.5,4.5) node {$M$};
	\draw [fill=eqeqeq] (-2.529065300318972,3.4716800159165246) circle (5.5pt);
	\draw[] (-2,4.2) node {$M_{13}$};
	\draw [fill=eqeqeq] (2.8584222402083728,-2.8992377899451656) circle (5.5pt);
	\draw[] (3.35013737287555,-2.1) node {$M_{34}$};
	\draw [fill=eqeqeq] (6.53552936824566,7.779021661238361) circle (5.5pt);
	\draw[] (7.6,7.5) node {$M_{23}$};
	\draw [fill=black] (0.928896264065991,2.088938128623705) circle (5.5pt);
	\draw[color=black] (0,3.1) node {$V_1$};
	\draw [fill=black] (2.1205227544153047,0.36162633527332094) circle (5.5pt);
	\draw[color=black] (3,.4) node {$V_2$};
	\draw [fill=black] (11.15588945987388,0.48777283504246294) circle (5.5pt);
	\draw[color=black] (10.5,-.2) node {$V_3$};
	\draw [fill=black] (0.43336099838949727,5.728559217902819) circle (5.5pt);
	\draw[color=black] (1.1694876540906722,6.) node {$V_4$};
	\draw [fill=black] (-5.246680832412158,0.9253023676946988) circle (5.5pt);
	\draw[color=black] (-6.2,1.4) node {$V_{123}$};
	\draw [fill=black] (-0.19131801582759858,2.3698108968846037) circle (5.5pt);
	\draw[color=black] (0.6,1.2) node {$V_{134}$};
	\draw[color=black] (-8.3,6) node {$c_{13}^4$};
	\draw[color=black] (-7.8,-0.8468615840299902) node {$c_{13}^2$};
	\draw[] (14.5,3) node {$q_3$};
	\draw[] (-7.8,3.2) node {$q_{123}$};
	\end{large}
	\end{tikzpicture}
}
\numberwithin{equation}{section}

\newtheorem{theorem}{Theorem}[section]
\newtheorem{lemma}[theorem]{Lemma}
\newtheorem{definition}[theorem]{Definition}

\begin{document}

\title{Miquel dynamics, Clifford lattices and the Dimer model}

\author{Niklas C. Affolter}
\address{Niklas Affolter, TU Berlin, Institute of Mathematics, Strasse des 17. Juni 136,
	10623 Berlin, Germany}
\email{affolter@math.tu-berlin.de}

\begin{abstract}{
		Miquel dynamics were introduced by Ramassamy as a discrete time evolution of square grid circle patterns on the torus. In each time step every second circle in the pattern is replaced with a new one by employing Miquel's six circle theorem. Inspired by these dynamics we define the Miquel move, which changes the combinatorics and geometry of a circle pattern locally. We prove that the circle centers under Miquel dynamics are Clifford lattices, considered as an integrable system by Konopelchenko and Schief. Clifford lattices have the combinatorics of an octahedral lattice and every octahedron contains six intersection points of Clifford's four circle configuration. The Clifford move replaces one of these circle intersection points with the opposite one.\\% We show that the set of integrable circle patterns as defined by Bobenko, Mercat and Suris is closed under both the Miquel and the Clifford move.\\
		We establish a new connection between circle patterns and the dimer model: If the distances between circle centers are interpreted as edge weights, the Miquel move preserves probabilities in the sense of of urban renewal.} 
\end{abstract}

\maketitle

% !TEX root = oupau.tex

\section{Introduction} \label{sec:intro}
Miquel dynamics were first introduced by Ramassamy following an idea of Kenyon, see \cite{miquel} and references therein. In each time step, this discrete dynamical system replaces every second circle of a square grid circle pattern. If the circle pattern is doubly periodic, it is conjectured that these dynamics feature a form of discrete integrability, and that they are related to dimer statistics or dimer integrable systems \cite{dimerintegrable}. First progress towards integrability has been made by Glutsyuk and Ramassamy in \cite{miquel2} for the case of the doubly periodic $2\times2$ grid.\\
In Theorem \ref{th:miquelclifford} we show that the collection of circle centers under Miquel dynamics form a special case of Clifford lattices, a discrete integrable system studied by Konopelchenko and Schief \cite{wolfgangclifford}. We introduce the geometric star-ratio function, and the centers of circle patterns are exactly the Clifford lattices with real star-ratios.\\
If the star-ratios are real and positive, we say a circle pattern is Kasteleyn. To a Kasteleyn circle pattern we associate a dimer model with edge weights equal to the distances of circle centers. We show in Theorem \ref{th:miqueldimers} that the Miquel move induces urban renewal on the associated dimer model. This proves that under Miquel dynamics the star-ratios transform exactly as the face weight coordinates of the discrete cluster integrable system as introduced by Goncharov and Kenyon \cite{dimerintegrable}. However, the star-ratios do not fully govern Miquel dynamics as they do not uniquely determine the circle pattern. The analysis of the additional information needed to fully describe Miquel dynamics is a topic of ongoing research.\\
We want to stress that all our proofs are of a local nature. This has the advantage that the theory is not restricted to $\Z^2$ combinatorics. %Thus, we show that integrable circle patterns in the sense of Bobenko, Mercat and Suris \cite{nonlinearconformal} are a closed reduction of Clifford lattices with respect to their intersection points and a closed reduction of Miquel dynamics with respect to their circle centers.\\

The structure of the paper is as follows: In the following section~\ref{sec:prelim} we revisit the definition of dimer statistics and introduce notation for the graphs, circle patterns and the star-ratios considered throughout this paper. We also define both the local Miquel and the local Clifford move. After stating the two main theorems \ref{th:miquelclifford} and \ref{th:miqueldimers} in section~\ref{sec:main} we introduce in section \ref{sec:mobius} the M\"obius geometric mutation map which underlies both the Clifford and the Miquel move. We also calculate the transformation formulas for star-ratios under the M\"obius mutation map. In section \ref{sec:clifford} we study the Clifford configuration, giving a geometric construction of the M\"obius mutation map and relating it to work of Konopelchenko and Schief \cite{wolfgangclifford}. Additionally, we derive several useful lemmas by investigating the relation between the Clifford configuration, integrable cross ratio systems and integrable circle patterns as defined by Bobenko, Mercat and Suris \cite{nonlinearconformal}. Finally, we assemble the pieces in section~\ref{sec:miquel} to relate star-ratios with the Miquel move thereby proving the two main theorems. Even though we prove everything locally, it is interesting how this translates to the setting of lattice dynamics, which we briefly outline in section~\ref{sec:lattice}. Section~\ref{sec:conclusion} contains concluding remarks and open questions.
% !TEX root = oupau.tex
\section{Preliminaries} \label{sec:prelim}

%If $P,P'\in \C$ are two different point then we will denote by $PP'$ the euclidean straight line through $P$ and $P'$. Additionally, if we work with a labeling of objects by indices 1 to $n$ then we consider $n+1$ equivalent to $1$ and $n+2$ equivalent to 2 et cetera.

\subsection{Dimer statistics}\label{subsec:dimers}

Given a graph $G = (V,E)$ we call a function $\omega: E\rightarrow \R_+$ an \emph{edge weight function}, where $\R_+$ is the set of strictly positive reals. Equivalently, we write $\omega\in \R_+^{E}$ and call $\omega$ the \emph{edge weights}.

\begin{definition} 
	A \emph{perfect matching} of a simple graph $G=(V,E)$ is a subset $M\subset E$ is a subset of the edge set such that each vertex of the graph is incident to exactly one edge in $M$. We denote the set of perfect matchings of a graph by $\mathcal M(G)$. The \emph{weight} $\omega(M)$ of a perfect matching $M$ with respect to edge weights $\omega$ is:
	\begin{align}
	\omega(M) = \prod_{e\in M} \omega(e)
	\end{align}
	The \emph{dimer partition function} $Z_G$ is defined :
	\begin{align}
	Z_G: \R_+^{E} \rightarrow \R_+,\qquad \omega \mapsto \sum_{M\in \mathcal M(G)}\omega(M)
	\end{align}
	The \emph{probability} $P_\omega(M)$ of a perfect matching is proportional to its weight with the partition function as normalization constant.
	\begin{align}
		P_\omega(M) = \f{\omega(M)}{Z_G(\omega)}
	\end{align}
\end{definition}

\subsection{Bipartite surface graphs $\obs$}
We start by defining the class $\obs$ of bipartite surface graphs for which we prove our lemmas and theorems.
\begin{definition} \label{def:surface}
	A graph $G = (V,E)$ embedded in an oriented and closed surface is in $\obs$ if:
	\begin{itemize}
		\item The complement of $G$ is a collection $F$ of disjoint open disks, and $(V,E,F)$ is a locally finite CW-decomposition of the surface.
		\item $G$ is bipartite, that is the set of vertices $V$ is the disjoint union of the two independent sets $V^+$ and $V^-$.
		\item Every vertex in $V$ has degree at least 3 and every face has degree at least 2.
	\end{itemize}
\end{definition}

\begin{figure}[th]
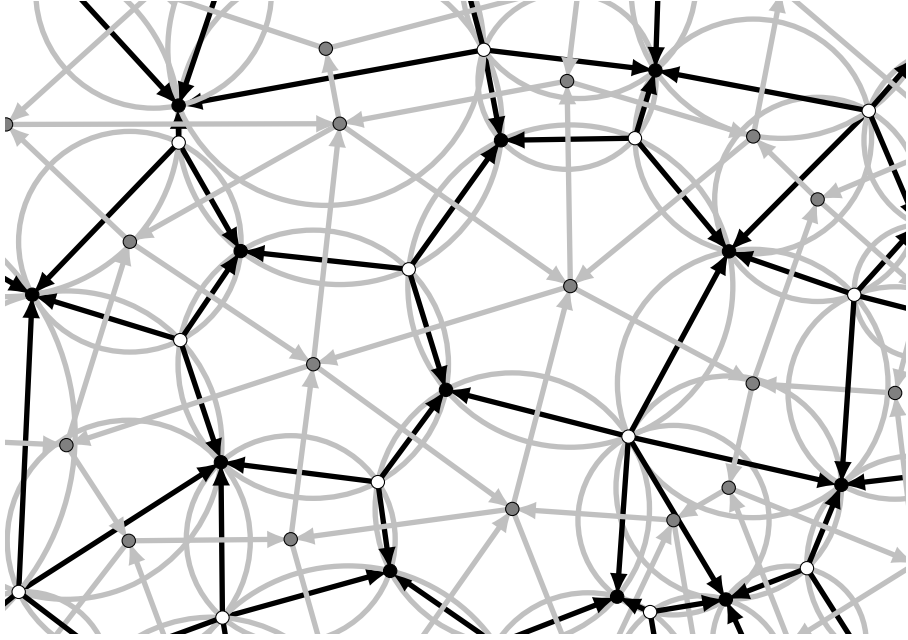
	
	\centering
	\ggborientationexample
	\caption{An excerpt of a circle pattern $z\in \CP^G$, $G\in \obs$. The vertices in $V^+$ ($V^-$) are colored black (white), the circle centers which correspond to faces of $G$ are colored gray. Edges of $G$ are drawn with black arrows and edges of $G^*$ with gray arrows. The direction of an arrow indicates the orientation of the edge.}
	\label{fig:orientationexample}
\end{figure}

In order to simplify notation we assume that we can identify an edge $e\in E$ with its two incident vertices $v$ and $v'$ and therefore write $e =(v,v')$. Similarly for dual edges $e^*\in E^*$ we write $e^* = (f,f')$ where $f$ and $f'$ are the two faces incident to $e$. We usually consider edges at some distinct face $f_0$, and the cyclic order of edges around $f_0$ will usually clear up any confusion.\\
We also use a standard orientation of the edges of both $G$ and the dual $G^*$. An edge $e=(v^-,v^+)$ is always incident to a vertex $v^+\in V^+$ and a vertex $v^-\in V^-$, and we orient that edge as pointing from $v^-$ to $v^+$. Each dual edge $e^*$ is oriented such that it crosses the oriented primal edge $e$ from the left to the right. As a consequence, the dual edges are oriented counter clockwise around vertices in $V^+$ and clockwise around vertices in $V^-$. For an example of the standard orientation see figure \ref{fig:orientationexample}.

\begin{figure}[th]
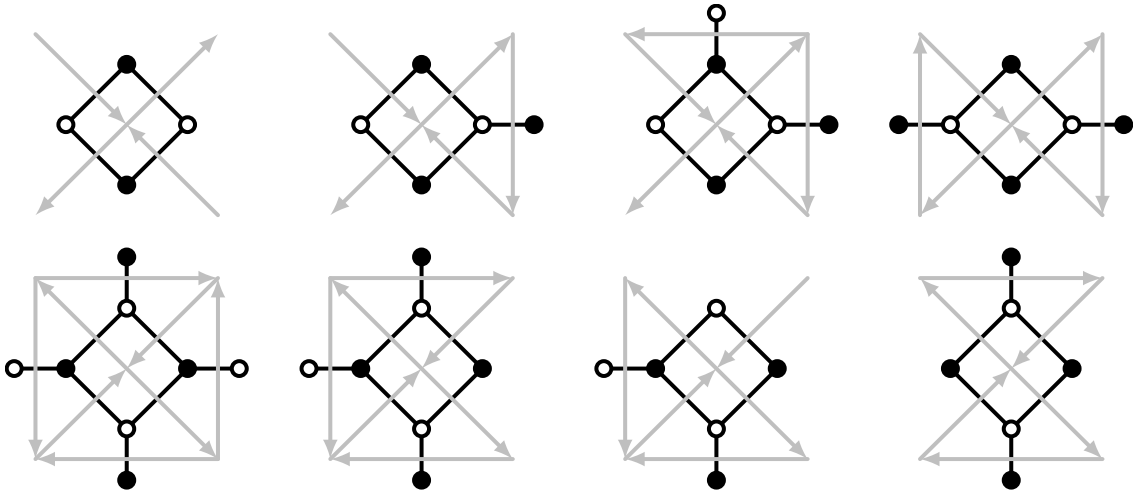
	
	\centering
	\begin{minipage}[b]{0.22\textwidth}
		\ggbmutationcomb00000
	\end{minipage}
	\begin{minipage}[b]{0.22\textwidth}
		\ggbmutationcomb10000
	\end{minipage}
	\begin{minipage}[b]{0.22\textwidth}
		\ggbmutationcomb11000
	\end{minipage}
	\begin{minipage}[b]{0.22\textwidth}
		\ggbmutationcomb10100
	\end{minipage}\\
	\begin{minipage}[b]{0.22\textwidth}
		\ggbmutationcomb11111
	\end{minipage}
	\begin{minipage}[b]{0.22\textwidth}
		\ggbmutationcomb01111
	\end{minipage}
	\begin{minipage}[b]{0.22\textwidth}
		\ggbmutationcomb00111
	\end{minipage}
	\begin{minipage}[b]{0.22\textwidth}
		\ggbmutationcomb01011
	\end{minipage}
	\caption{The eight local configurations up to rotations at a face $f$ with four neighbours. The black edges are the edges in $N_f$ of the primal graph $G\in \obs$. The gray arrows are the edges of the dual graph $G^*$. In each column the bottom configuration is the mutation of the upper one and vice versa. Under mutation the set of faces is preserved, and the outer four faces are connected to the rest of $G$.}
	\label{fig:mutationcombinatorics}
\end{figure}

\begin{definition} \label{def:edgeneighbours}
	Let $G\in \obs$ and let $f \in F$ be a quadrilateral with the four neighbours $f_1,f_2,f_3,f_4$. Then we define the \emph{edge neighbourhood $N_f$} of $f$ as the set of edges that is not incident with any other face except $f_1,f_2,f_3,f_4$ or $f$. The dual edge neighbourhood $N^*_f$ is the set of edges dual to edges in $N_f$.
\end{definition}

\begin{definition}
	Given $G\in\obs$ and a quadrilateral $f\in F$ we denote by $\mut_fG\in \obs$ the graph resulting from the \emph{4-mutation at $f$}. It differs from $G$ by a local change of combinatorics centered at $f$. We will write $\tilde G = \mut_f G$ and $\tilde G = (\tilde V,\tilde E, \tilde F)$. The face set is invariant under mutation $\tilde F = F$, and therefore it is easiest to describe the change in combinatorics with respect to the dual of $G$:
	\begin{align}
		\tilde E^* = E^*\ \Delta\ \menge{(f_1,f_2),(f_3,f_2),(f_3,f_4),(f_4,f_1)}
	\end{align}
	Where $\Delta$ denotes the symmetric difference operator. Figure \ref{fig:mutationcombinatorics} shows all eight possible configurations and their relation by mutation. The mutation at $f$ only affect edges in $N_f$ and only inserts or deletes vertices incident to $f$. The mutation is an involution i.e., $G = \mut_f \tilde G$.
\end{definition}

The term mutation is borrowed from the theory of cluster algebras, where the 4-mutation corresponds to the mutation at a degree four vertex. Indeed, both dimer statistics and Miquel dynamics may be formulated in the language of cluster algebras. We refrain from doing so because we do not need to use any results from the theory of cluster algebras in this paper.

\subsection{Urban renewal}

Because mutation at a quadrilateral $f$ only affects the edges in $N_f$, the complementary sets $N^\mathsf{c}$ and $\tilde N^\mathsf{c}$ can be identified. Therefore we are able to compare a matching $M_0\in \mathcal M(G)$ with a matching $M$ from $\mathcal M(G)$ or $\mathcal M(\tilde G)$ by comparing their restrictions to $ N_f^\mathsf{c}$. If $M_0$ and $M$ agree on $N_f^\mathsf{c}$, we write $M_0 \sim_f M$.
%\begin{align}
%M \sim_f M' \qquad \eq \qquad M\cap N_f^\mathsf{c} =  M'\cap  
%\end{align}

\begin{definition} \label{def:urbanrenewal} Fix a quadrilateral $f\in F$ and consider two edge weight functions $\omega \in \R_+^{E}$ and  $\tilde\omega \in \R_+^{\tilde  E}$. We say they are \emph{related by urban renewal} if the following two conditions are satisfied:\\
\begin{itemize}
	\item[(i)] For all $e\in N_f^\mathsf{c}: \quad \omega(e) = \tilde\omega(e)$ 
	\item[(ii)] For any fixed matching $M_0\in \mathcal{M}(G)$:
	\begin{align}
		\sum_{\substack{M \in\mathcal{M}(G)\\M \sim_f M_0 } }P_\omega(M) = \sum_{\substack{M \in\mathcal{M}(\tilde G)\\M \sim_f M_0 } }P_{\tilde\omega}(M)
	\end{align}
	
\end{itemize}
\end{definition}

\begin{definition}
	For a graph $G\in \obs$ we define the \emph{face weight function}  $\tau$ as an alternating ratio of the edge weights as follows:
	\begin{align}
	\tau: \R_+^{E} \rightarrow \R_+^{F}, \qquad (\tau(\omega))(f) = \prod_{\substack{e\in E\\e^*=(f,f')}} \omega(e) \prod_{\substack{e\in E\\e^*=(f',f)}} \left(\omega(e)\right)^{-1}
	\end{align}	
	The first product is over the edges whose duals point into $f$ and the second product is over edges  whose duals point away from $f$, where we use the orientation of the dual edges described in definition \ref{def:surface}.
\end{definition}

The next lemma relating urban renewal and face weights is well known \cite{propp, dimerintegrable}.
\begin{lemma} \label{lem:urbanfaceweights}
	Consider two edge weight functions $\omega \in \R_+^{E}$ and $\tilde\omega \in \R_+^{\tilde E}$ that agree on all $e\in N^\mathsf{c}$. Adopt the following abbreviations: $\tau_f = (\tau(\omega))(f)$ and $\tilde\tau_f = (\tau(\tilde\omega))(f)$. Then $\omega$ and $\tilde{\omega}$ are related by urban renewal if and only if the face weights $\tau$ before and the face weights $\tilde \tau$ after mutation are related as follows:
	\begin{align}
		\tilde\tau_{f'} = \begin{cases}
			\tau^{-1}_{f'} & f'=f\\
			\tau_{f'} (1+ \tau_f) & (f',f) \in E_G\\
			\tau_{f'} (1+ \tau^{-1}_f)^{-1} & (f, f') \in E_G\\
			\tau_{f'} & \mbox{else}\\
		\end{cases}
	\end{align}
\end{lemma}

\subsection{Circle patterns}

\begin{definition}
	Let $G$ be a simple graph and $z: V_G\rightarrow \Ce$. Then we call $z$ a \emph{vertex drawing of $G$} if no two adjacent vertices are mapped to the same point. Similarly, if $z: F_G\rightarrow \Ce$ is a map such that no two adjacent faces are mapped to the same point we call this a \emph{face drawing of $G$}.
\end{definition}

We work with circles in $\Ce$, which encompass both the euclidean circles and the straight lines of $\R^2$ which we identified with $\C = \Ce \setminus \menge{\infty}$. We call the reflection of $\infty$ about a circle the circle center. Therefore the center of a straight line is $\infty$.

\begin{definition}
	Let $z$ be a vertex drawing of $G\in \obs$. We call $z$ a \emph{circle pattern} if two conditions are fulfilled: The first condition is that for each face $f\in F$, the vertices incident to $f$ are mapped to a common circle in $\Ce$. This allows us to define a map $z^* \in \Ce^F$ that maps each face of $G$ to the corresponding circle center in $\Ce$. The second condition for $z$ to be a circle pattern is that $z^*$ is a face drawing of $G^*$. We denote the set of circle patterns of $G$ by $\CP^G$.
\end{definition}

For an example of a circle pattern see figure \ref{fig:orientationexample}.

\begin{definition} 
Let $G\in \obs$, then the \emph{star-ratio} $\sr$ is a map from $\Ce^F$ to itself:
%$Q\in \obs^*$, then the \emph{star-ratio} $\sr$ is defined as a map between function on the vertices of $V$.
\begin{align}
	\sr: \Ce^F \rightarrow \Ce^F,\quad (\sr(z))(v) = -\f{\prod_{(v', v) \in E} (z(v')-z(v))}{\prod_{(v, v') \in E} (z(v')-z(v))}
\end{align}
For computations we will also sometimes write:
\begin{align}
	\sr(y;y_1,y_2,y_3,y_4) = -\f{(y_1-y)(y_3-y)}{(y_2-y)(y_4-y)}
\end{align}
Where $y,y_1,y_2,y_3,y_4$ are all in $\Ce$.
\end{definition}
In this definition we use the standard approach of M\"obius geometry in which infinities cancel and division by zero yields infinity.

\begin{lemma} \label{lem:cphaverealsr}
Given $G\in \obs$ and a circle pattern $z\in \CP^G$, the star-ratios of the dual drawing $z^*\in \Ce^F$ are real. That is:
\begin{align}
	\forall f\in F: \qquad (\sr(z^*))(f) \in \R
\end{align}
\end{lemma}
\proof{Let $v_1, v_2,...,v_{2n}$ be the vertices incident to $f$ and $f_1, f_2,...,f_{2n}$ the incident faces, then $z(v_{k+1})$ is the reflection of $z(v_k)$ about the line through $z^*(f)$ and $z^*(f_k)$. Because $v_{2n+1} = v_{1}$, the concatenation of all reflections is the identity. The argument of the star-ratio is exactly $\pi$ plus half of the sum of the reflection angles and is therefore in $\pi \Z$.\qed\\	
}
\begin{definition}
	We call a circle pattern $z\in \CP^G$ of $G\in\obs$ \emph{Kasteleyn} if and only if all the star-ratios of $z^*\in \Ce^F$ are real positive.
\end{definition}
The usage of the term Kasteleyn is justified because if one would use the complex vectors of the oriented differences between circle centers as weights in an adjacency matrix, this matrix would be Kasteleyn-flat \cite{kasteleyn,kuperberg}. The square root of the determinant of a Kasteleyn-matrix counts weighted perfect matchings.\\ Furthermore, any circle pattern $z\in \CP^G$ such that the drawing $z$ is an embedding of $G$ into $\C$ is Kasteleyn. Therefore a natural and well investigated class of Kasteleyn circle patterns does indeed exist.\\

\begin{figure}[h]
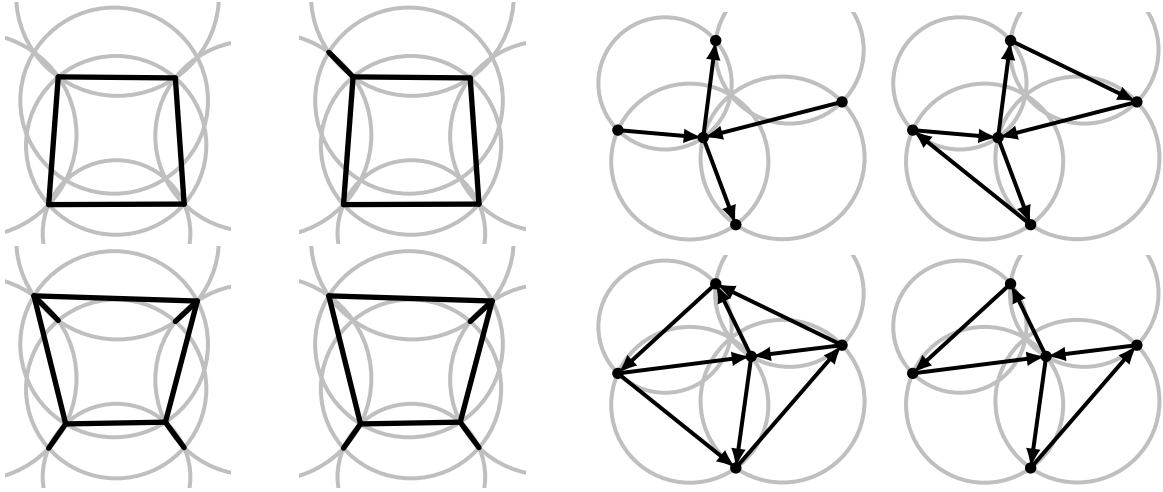
	
	\centering
	\begin{minipage}[b]{0.22\textwidth}
		\ggbmiqcomb00000
	\end{minipage}
	\begin{minipage}[b]{0.22\textwidth}
		\ggbmiqcomb10000
	\end{minipage}
%	\begin{minipage}[b]{0.22\textwidth}
%		\ggbmiqcomb11000
%	\end{minipage}
%	\begin{minipage}[b]{0.22\textwidth}
%		\ggbmiqcomb10100
%	\end{minipage}\\
		\begin{minipage}[b]{0.22\textwidth}		
			\ggbcliffordcomb00000
		\end{minipage}
		\begin{minipage}[b]{0.22\textwidth}		
			\ggbcliffordcomb10100
		\end{minipage}\\
	\begin{minipage}[b]{0.22\textwidth}
		\ggbmiqcomb11111
	\end{minipage}
	\begin{minipage}[b]{0.22\textwidth}
		\ggbmiqcomb01111
	\end{minipage}
%	\begin{minipage}[b]{0.22\textwidth}
%		\ggbmiqcomb00111
%	\end{minipage}
%	\begin{minipage}[b]{0.22\textwidth}
%		\ggbmiqcomb01011
%	\end{minipage}\\
		\begin{minipage}[b]{0.22\textwidth}
			\ggbcliffordcomb11111
		\end{minipage}
		\begin{minipage}[b]{0.22\textwidth}
			\ggbcliffordcomb01011
		\end{minipage}
	\caption{On the left we have two examples of configurations before (on top) the Miquel move and after (below) the Miquel move. The edges of $G$ are drawn in black, notice that the quadrilaterals have a different set of legs. On the right side we see two examples of the Clifford move, where the arrows indicate the oriented dual edges.}
	\label{fig:geometrymoves}
\end{figure}

\begin{definition}
	Let $G\in \obs$ and $z$ be a face drawing of $G$. Then we call a face $f\in F$ a \emph{valid face} if the following two conditions are fulfilled:
	\begin{itemize}
		\item The face $f$ has exactly four neighbours $f_1,f_2,f_3,f_4$, which are listed here in the cyclic order with respect to the orientation of the graph.
		\item No two consecutive faces are mapped to the same point, that is $z(f_k) \neq z(f_{k+1})$ for all four values of $k$.
	\end{itemize}
\end{definition}

We will now introduce the Miquel and the Clifford move which are based on Miquel's six circle theorem respectively Clifford's four circle theorem. They locally alter the geometry and are accompanied by a mutation of the combinatorics. For examples see figure \ref{fig:geometrymoves}.
Miquel's six circle theorem states that given the four circles $c_1,c_2,c_3,c_4$ and their four pairs of intersection points $\menge{I_k,\tilde I_k} = c_k\cap c_{k+1}$, then: If the four points $I_1,I_2,I_3,I_4$ are on a circle $c$  then $\tilde I_1,\tilde I_2,\tilde I_3,\tilde I_4$ are on a circle $\tilde c$. The Miquel move replaces the circle $c$ with the circle $\tilde c$. For the definition of the Miquel move, denote the center of $c$ by $M$ and the center of $\tilde c$ by $\tilde M$.
\begin{definition}\label{def:miquelmove} 
Let $G\in \obs$, $z\in \CP^G$ be a circle pattern and $f\in F$ a valid face with the four neighbours $f_1,f_2,f_3,f_4$, such that $z$ does not map the five points $f,f_1,f_2,f_3,f_3$ to a common line. Let $v_1,v_2,v_3,v_4$ be the four vertices incident to $f$. Identify $z^*(f) = M$  with the center of circle $c$ in Miquel's six circle theorem. Also identify $z(v_k)$ with $I_k$. Then the \emph{Miquel move} replaces the circle associated to $f$ with the alternate one that exists due to Miquel's theorem:
\begin{align}
	\miq_f&: \CP^G \rightarrow \CP^{\tilde G}\\
	 \qquad (\miq_f(z^*))(f') &= \begin{cases}
	\tilde M & f' =f\\
	z^*(f') & f' \neq f
	\end{cases}\\
	\qquad (\miq_f(z))(v') &= \begin{cases}
	\tilde I_k & v' = v_k\\
	z(v') & v' \neq v
	\end{cases}
\end{align}
\end{definition}

Clifford's four circle theorem considers circles $c_1$ to $c_4$ that intersect in one common point $I$. Therefore the four pairs of intersection points of subsequent circles are $\menge{I,J_{k,k+1}} = c_k\cap c_{k+1}$. There are two more intersection points $J_{k,k+2} = c_k\cap c_{k+2}$. Define the additional four circles $\tilde  c_k = (J_{k,k+1},J_{k-1,k}, J_{k-1,k+1})$. Clifford's four circle theorem states that in fact all four circles $\tilde c_k$ intersect in one point $\tilde  I$.
\begin{definition}\label{def:cliffordmove}  Let $G \in \obs$ be a graph and $z\in \Ce^F$ be a face drawing of $G$ and $f\in F$ a valid face with the four neighbours $f_1,f_2,f_3,f_4$, such that $z$ does not map the five points $f,f_1,f_2,f_3,f_3$ to a common circle. Identify $z(f) = I$ and $z(f_k) = J_{k,k+1}$ with the intersection points in Clifford's four circle theorem above, then set $(\cli_f(z))(f) = \tilde I$. This defines the \emph{Clifford move}.
\begin{align}
	\cli_f: \Ce^F \rightarrow \Ce^{\tilde F}: \qquad (\cli_f(z))(f') = \begin{cases}
		\tilde I & f' =f\\
		z(f') & f' \neq f
	\end{cases}
\end{align}
\end{definition}

% !TEX root = oupau.tex

\section{Main theorems} \label{sec:main}

\begin{theorem} \label{th:miquelclifford}
	Let $G \in \obs$ be a graph, let $z\in \CP^G$ be a circle pattern of $G$ and let $z^* \in \Ce^F$ be the dual drawing induced by the circle centers of $z$. Let $f\in F$ be a valid face with the four neighbours $f_1,f_2,f_3,f_4$, such that $z^*$ does not map the five points $f,f_1,f_2,f_3,f_3$ to a common circle. Then the Miquel move at $f$ acts on the circle pattern such that the circle centers change as they do under the Clifford move:
	\begin{align}
		(\miq_f(z))^* = \cli_{f}(z^*)
	\end{align}

\end{theorem}
It is an implication of this theorem that the change of the circle centers does not depend on a particular choice of circles with these circle centers. That is, if $z_1,z_2 \in \CP^G$ are such that $z_1^*=z_2^*$, then:
\begin{align}
	(\miq_f(z_1))^* = (\miq_f(z_2))^*
\end{align}
	
\begin{definition}\label{def:geometricdimer}
Given a graph $G \in \obs$, the map $\psi:\CP^G\rightarrow\R_+^E$ associates edge weights to a circle pattern $z\in\CP^G$ as follows.
\begin{align}
	\psi:\CP^G\rightarrow\R_+^E,\qquad (\psi(z))(e) = |z^*(f)-z^*(f')|\quad \forall e\in E,\ (f,f')=e^*
\end{align}
\end{definition}

\begin{theorem} \label{th:miqueldimers}
	Let $G \in \obs$ be a surface graph, $z\in \CP^G$ a Kasteleyn circle pattern and $f \in F$ a valid face  with the four neighbours $f_1,f_2,f_3,f_4$, such that $z^*$ does not map the five points $f,f_1,f_2,f_3,f_3$ to a common line. Then the edge weights $\omega = \psi(z) \in \R_+^E$ and $\tilde\omega = \psi(\miq_fz) \in \R_+^{\tilde E}$ are related by urban renewal. %That is, the probability measures of the associated dimer models on $G$ and $\tilde G$ fulfill for any fixed matching $M_0\in \mathcal{M}(G)$:
	%\begin{align}
	%	\sum_{\substack{M \sim_f M_0 \\M \in\mathcal{M}(G)} }P_\omega(M) = \sum_{\substack{M \sim_f M_0 \\M \in\mathcal{M}(\tilde G)} }P_{\tilde\omega}(M)
	%\end{align}
	%The sums are over all matchings that agree with $M_0$ on $N^\mathsf{c}$, i.e. the edge set that is invariant under the mutation as explained in subsection \ref{subsec:dimers}.
\end{theorem}

We prove Theorem \ref{th:miquelclifford} and Theorem \ref{th:miqueldimers} at the end of section \ref{sec:miquel}.
% !TEX root = oupau.tex

\section{The star-ratio preserving M\"obius map} \label{sec:mobius}
In this section we introduce M\"obius maps that preserve star-ratios. In particular we are interested in the following problem: Given five points $z,z_1,z_2,z_3,z_4 \in \Ce$, where are all points $\tilde z \in \Ce$ such that the star-ratio of $z,z_1,z_2,z_3,z_4$ and $\tilde z,z_1,z_2,z_3,z_4$ are the same? It turns out that there is at most one such $\tilde z \neq z$, and that $\tilde z$ can be expressed as a fractional linear transform of $z$ with coefficients consisting of $z_1,z_2,z_3$ and $z_4$. This motivates the definition a unique non-trivial M\"obius map that takes any $z$ to $\tilde z$ for fixed $z_1,z_2,z_3,z_4$.

\begin{definition} 
Given four points $z_1,z_2,z_3,z_4$ such that $z_k\neq z_{k+1}$, we define the \emph{M\"obius mutation map} $\mob(z_1,z_2,z_3,z_4)$. 
\begin{align}
	\mob(z_1,z_2,z_3,z_4): \Ce \rightarrow \Ce,\qquad z \mapsto\f{zC_2 +C_3}{zC_1-C_2 }
\end{align}
The coefficients are the following homogeneous polynomials in the four points $z_k$:
\begin{align}
	C_1 &= z_1-z_2+z_3-z_4\\
	C_2 &= z_1z_3-z_2z_4\\
	C_3 &= z_2z_3z_4-z_1z_3z_4+z_1z_2z_4-z_1z_2z_3 %z_1z_2z_3z_4\left(z_1^{-1}-z_2^{-1}+z_3^{-1}-z_4^{-1}\right)
\end{align}
\end{definition}
Of course, a fractional linear transform only corresponds to a M\"obius map if the determinant given by ${-C_2^2-C_1C_3}$ is not zero. However, the determinant is non vanishing because it factors as follows:
\begin{align}
	\det \mob(z_1,z_2,z_3,z_4) = (z_1-z_2)(z_2-z_3)(z_3-z_4)(z_4-z_1)
\end{align}
Notice that the trace of the M\"obius mutation map is zero, consequently it can not be the identity.
If $z_1\neq z_3$ or $z_2 \neq z_4$ then the M\"obius mutation map is the unique map that maps $z_k$ to $z_{k+2}$. It is straight-forward to verify that the map interchanges the points as stated. Because we have prescribed the image of at least three points, the map is also unique. If $z_1 = z_3$ and $z_2 = z_4$ then the mutation map is the unique map that has $z_1$ and $z_2$ as fixed points and that sends the point $\f12(z_1+z_2)$ to the point at infinity.

\begin{lemma} \label{lem:srmobius}
	Let $z_1,z_2,z_3,z_4 \in \Ce$ such that $z_k\neq z_{k+1}$ and let $z\in \Ce$ be such that $z\neq z_k$ for any of the four points. Then the M\"obius mutation map and the identity are the unique M\"obius maps $M: \Ce\rightarrow\Ce$ such that the following relation holds:
	\begin{align}
		\f{\sr(z;z_1,z_2,z_3,z_4)}{\sr(M(z);z_1,z_2,z_3,z_4)}=1\label{eq:srsame}
	\end{align}
\end{lemma}
\proof{Notice that the quotient of star-ratios is decomposable into cross ratios as follows:
	\begin{align}
		\f{\sr(z;z_1,z_2,z_3,z_4)}{\sr(M(z);z_1,z_2,z_3,z_4)} =\cro(z,z_1,M(z),z_2)\cro(z,z_3,M(z),z_4)
	\end{align}
	Therefore, if $M$ satisfies equation  \eqref{eq:srsame}, it will do so also after some other M\"obius map has been applied. Hence we can assume that no point is at infinity. In that case, equation \eqref{eq:srsame} is a quadratic equation for $M(z)$ and a simple calculation shows that indeed $M = \mob(z_1,z_2,z_3,z_4)$ or $M$ is the identity. These are the only possibilities because there are at most two solutions.\qed
}

\begin{definition} Let $G\in \obs$, $z$ a face drawing and let $f\in F$ be a valid face of $G$. Extend the definition of the M\"obius mutation map to the setting of surface graphs and mutations as follows:
\begin{align}
	\mob_f: \Ce^F \rightarrow \Ce^{\tilde F}, \qquad(\mob_f(z))(f') = \begin{cases}
		\mob(z(f_1),z(f_2),z(f_3),z(f_4))(z(f)) & f' = f\\
		z(f') & f' \neq f
	\end{cases}
\end{align}

\end{definition}

\begin{figure}[ht]
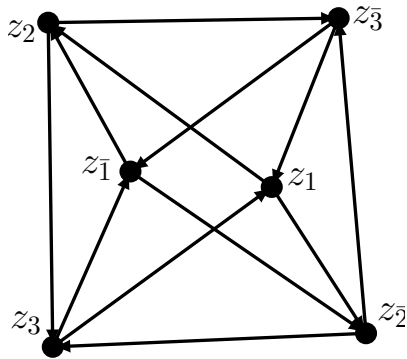

	\centering
	\ggbsroctahedron
	\caption{The octahedron which is reflected onto itself by the M\"obius mutation map $\mob(z_1,z_2,z_{\bar{1}}, z_{\bar{2}})$.}
	\label{fig:octahedronquiver}
\end{figure}

\begin{lemma} \label{lem:srmutation} Let $G\in \obs$, $z$ a face drawing and let $f\in F$ be a valid face. Let the tilde denote the quantities on the mutated graph $\widetilde G = \mut_f G$. Under M\"obius mutation, the star-ratios change according to the following relations: 
	\begin{align}
		\widetilde\sr(\tilde z)(f') = \begin{cases}
			(\sr(z))^{-1}(f') & f'=f\\
			\sr(z)(f') (1+ \sr(z)(f)) & (f',f) \in E\\
			\sr(z)(f') (1+ (\sr(z)(f))^{-1})^{-1} & (f, f') \in E\\
			\sr(z)(f') & \mbox{else}\\
		\end{cases}
	\end{align}	
\end{lemma}

\proof{Let $f \in F$ be as in the lemma with $(f_1,f),(f,f_2),(f_3,f),(f,f_4) \in E^*$. Introduce the 3-cube $H_3 \in \obs$, its dual is the octahedron. The faces of $H_3$ are $F_O = \menge{g_1,g_{\bar1},g_2,g_{\bar2},g_3,g_{\bar3}}$ as in figure \ref{fig:octahedronquiver}. Choose a face drawing $z'\in \Ce^{F_O}$ by writing $z'_k=z(u_k)$ and by imposing the following relations: 
\begin{align}
	 z'_2 &= z(f_1),\ z'_3 = z(f_2),\ z'_{\bar2} = z(f_3),\ z'_{\bar3} = z(f_4)\\
	z'_1 &= z(f),\ z'_{\bar1} = \mob(z'_2, z'_3,z'_{\bar2}, z'_{\bar3})(z'_1)
\end{align} 
By the definition of the M\"obius mutation map we can deduce the following symmetric action of the mutation map on the octahedron.
	\begin{align}
		z'_{\bar{k}}&=\mob(z'_l,z'_m,z'_{\bar{l}}, z'_{\bar{m}})(z'_k)\qquad \forall k,l,m \in \menge{1,2,3},\ l\neq m
	\end{align}
	Indeed, it does not matter which two pairs of opposite points we choose to define the mutation map.
	\begin{align}
		\mob(z'_k,z'_l,z'_{\bar{k}}, z'_{\bar{l}}) = \mob(z'_m,z'_n,z'_{\bar{m}}, z'_{\bar{n}})\qquad\forall k,l,m,n \in \menge{1,2,3},\ k\neq l,\ m\neq n
	\end{align}
	
	Denote by $\sr'_k = (\sr(z'))(g_k)$ the star-ratios of the cube and by $\sr_k = (\sr(z))(f_k)$ the star-ratios in the graph $G$.
	Lemma \ref{lem:srmobius} already ensures that $\sr'_k = \sr'_{\bar k}$. Additionally, a direct calculation shows that if we know $\sr_1$ we also know $\sr_2$ and $\sr_3$:
\begin{align}
	\sr'_3 = -(1+(\sr'_1)^{-1}),\quad \sr'_2 = -(1+\sr'_1)^{-1}\label{eq:sroctahedron}
\end{align}
	
	We identify the star-ratios $\sr$ of the cube with the star-ratios $\widetilde{\sr}$ before and after mutation at $f$ in $G$.
\begin{align}
	\sr_1 & = \sr'_1\\
	\widetilde\sr_1 &= (\sr'_{\bar1})^{-1}\\
	\f{\sr_1}{\widetilde\sr_1} = \f{\sr_3}{\widetilde\sr_3} &= {-\sr'_3}\\
	\f{\sr_2}{\widetilde\sr_2} = \f{\sr_4}{\widetilde\sr_4} &= {-\sr'_2}
	\end{align}	
	Inserting the expressions from equation \eqref{eq:sroctahedron} yields the statement of the lemma. \qed\\
}

Let us summarize what we have proven in this section. We have shown that for four points $z_1,z_2,z_3,z_4$ in valid position there exists a unique non trivial M\"obius map that preserves star-ratios with respect to the four points. If we use this map together with a mutation of combinatorics, we know that the star-ratios transform as in Lemma \ref{lem:srmutation}. \\It is our goal to identify star-ratios and face weights of dimer statistics. Due to the results of this section, it will suffice to show that the maps of Clifford and Miquel moves preserve the star-ratio of the face where the move is based at.
 
% !TEX root = oupau.tex
\section{Clifford configurations}\label{sec:clifford}
The combinatorial n-hypercube $H_n$ has as set of vertices all the subsets of $\menge{1,2,\dots,n}$. There is an edge between $I$ and $I'$ in the hypercube if these two sets differ by exactly one index. We call a vertex of the hypercube even (odd) if the cardinality of the corresponding index set is even (odd). The distinction between odd and even sets depends on what vertex has been labeled by the empty set. However, whether two vertices belong to the same parity is independent of that choice.

\begin{figure}[ht]
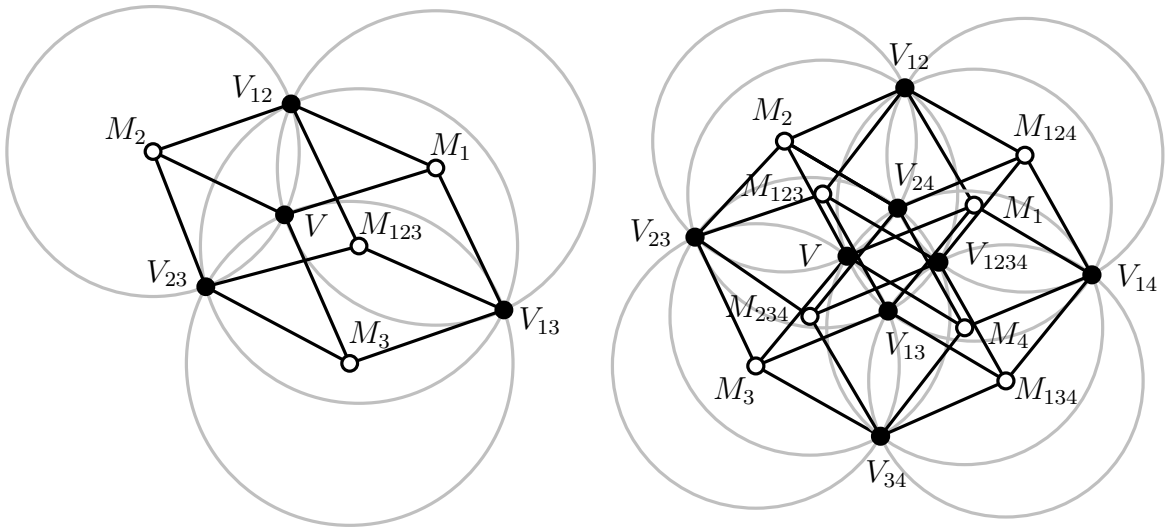

	\centering
	\begin{minipage}[b]{0.46\textwidth}
		\ggbcliffordthree
	\end{minipage}
	\begin{minipage}[b]{0.46\textwidth}
		\ggbcliffordfour
	\end{minipage}	
	\caption{The Clifford configurations $C_3$ and $C_4$.}
	\label{fig:clifford}
\end{figure}

\begin{figure}[h]
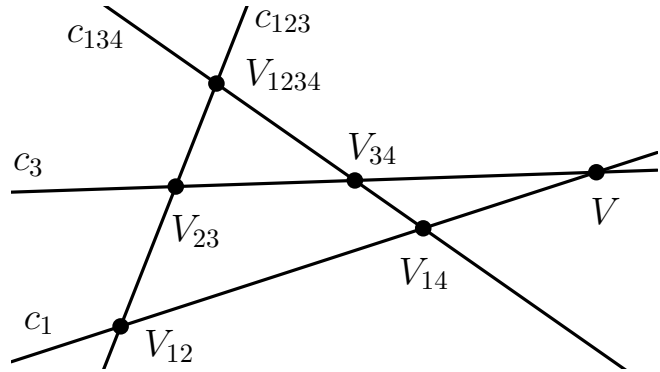

	\centering
	\ggbmenelaus
	\caption{Menelaus's theorem applied to $C_4$ with $V_{13}$ sent to infinity via a suitable M\"obius transformation.}
	\label{fig:menelaus}
\end{figure}

\begin{definition}
The \emph{Clifford-n-circle configuration} $C_n$ is a map from the n-hypercube $H_n$ to points and circles in the plane $\Ce$. Every even vertex is mapped to a point and every odd one to a circle, such that a point incident to a circle in the cube is mapped to that circle in the plane (see figure \ref{fig:clifford} for two examples).
\end{definition}
Clifford's theorem then says that given $n$ cyclically ordered circles in the plane all intersecting in one point, they can be extended uniquely to the whole Clifford configuration $C_n$.
We will denote circles in a Clifford configuration by $c_I$, intersection points by $V_I$ and circle centers by $M_I$. In order to facilitate notation, we will be slightly imprecise and identify the vertices of $H_n$ with the corresponding intersection points or circle centers in $\Ce$.
The connection of the Clifford configuration $C_4$ to the octahedron recurrence and the M\"obius mutation map becomes apparent in the next lemma.
\begin{lemma} \label{lem:cliffordmenelaus}
Consider the Clifford configuration $C_4$. The following two equivalent equalities hold:
\begin{align}
	\mob(V_{12},V_{23},V_{34},V_{14})(V) &= V_{1234}	\\
	\sr(V; V_{12},V_{23},V_{34},V_{14}) &= \sr(V_{1234}; V_{12},V_{23},V_{34},V_{14})
\end{align}

\end{lemma}
\proof{The proof follows the one of Konopelchenko and Schief \cite{wolfgangclifford} closely and employs Menelaus's theorem. It states that given three lines with three distinct intersection points $A_{12},A_{23},A_{13}$ and an additional line intersecting the three previous lines in the points $B_1,B_2,B_3$, we have an equality for the multi-ratio of the six points:
\begin{align}
	\mr(A_{12},B_2,A_{23},B_3,A_{13},B_1) := \f{A_{12}-B_2}{B_2-A_{23}}\f{A_{23}-B_3}{B_3-A_{13}}\f{A_{13}-B_1}{B_1-A_{12}} = -1
\end{align}

To apply this to $C_4$ we send $V_{13}$ to infinity via a M\"obius transformation. As a result the circles $c_1,c_3,c_{123}$ and $c_{134}$ are four lines with the six intersection points $V,V_{12},V_{14},V_{23},V_{34},V_{1234}$ and of course $V_{13}$ at infinity, compare with figure \ref{fig:menelaus}. In this setting Menelaus's theorem gives the following two equalities:
\begin{align}
	m_1&:=\mr(V_{34},V_{14},V_{1234},V_{12},V_{23},V)=-1\\
	m_2&:=\mr(V_{14},V_{34},V_{1234},V_{23},V_{12},V)=-1
\end{align}
Both $m_1$ and $m_2$ are actually M\"obius invariant quantities as they can be expressed in terms of cross ratios:
\begin{align}
	m_1 &= \cro(V_{34},V_{14},V_{1234},V_{13})\cro(V_{1234},V_{12},V_{23},V_{13})\cro(V_{23},V,V_{34},V_{13})\\
	m_2 &= \cro(V_{14},V_{34},V_{1234},V_{13})\cro(V_{1234},V_{23},V_{12},V_{13})\cro(V_{12},V,V_{14},V_{13})
\end{align}
The quotient of star-ratios on the other hand can be expressed by the quantities $m_1$ and $m_2$
\begin{align}
	\f{\sr(V; V_{12},V_{23},V_{34},V_{14})}{\sr(V_{1234}; V_{12},V_{23},V_{34},V_{14})}	= \f{m_2}{m_1} = 1
\end{align}
Because that expression is an equation in cross ratios it is invariant under M\"obius transformations and the claim is proven.\qed\\}

To simplify notation, let us introduce the shift operator $\shift J$ which acts as a symmetric difference operator on the indices of the arguments of a function, that is:
\begin{align}
	\shift J f(x_{I_1} , x_{I_2}, ..., x_{I_m}) &= f(x_{J \Delta I_1} , x_{J \Delta I_2}, ..., x_{J \Delta I_m})\\
	I\ \Delta\ J &= (I \cup J) \setminus (I\cap J)
\end{align}
Here, the $x$ variables could be vertices $V$, circles $c$ or circle centers $M$. Employing this notation the M\"obius mutation map acts as a shift operator on $C_4$.
\begin{align}
	\mob(V_{12},V_{23},V_{34},V_{14})(V_I) = \shift{1234} V_I,\qquad\mob(V_{12},V_{23},V_{34},V_{14})(c_I) = \shift{1234} c_{I}% \label{eq:M\"obiusshift}
\end{align}
This tells us that $\mob(V_{12},V_{23},V_{34},V_{14})$ acts as a total shift on the circles and intersection points of $C_4$. However, one should be aware of the fact that the M\"obius mutation will in general not act as a shift operator on the circle centers.\\
In terms of shift notation it is also possible to express the invariance of the star-ratio as follows:
\begin{align}
	\sr(V; V_{12},V_{23},V_{34},V_{14}) =  \shift{1234} \sr(V; V_{12},V_{23},V_{34},V_{14})
\end{align} 

Let us reiterate that in the $C_4$ configuration two points share the same star-ratio with respect to their same four neighbours. We can translate this fact into a lemma.

\begin{lemma} \label{lem:cliffordmobius}
	Let $G\in \obs$ be a graph, $z$ a face drawing of $G$ and let $f\in F$ be a valid face with the four neighbours $f_1,f_2,f_3,f_4$, such that $z$ does not map the five points $f,f_1,f_2,f_3,f_3$ to a common circle. The Clifford move (see definition \ref{def:cliffordmove}) coincides with the M\"obius mutation map.
	\begin{align}
		\cli_f(z) = \mob_f(z)
	\end{align}
\end{lemma}

This identification makes it natural to interpret the star-ratio preserving M\"obius mutation map as the extension of the Clifford move to the case where all five points lie on a common circle. In the next section we prove a similar result for the Miquel move, thereby proving Theorem  \ref{th:miquelclifford} relating Miquel dynamics and Clifford lattices.\\

%We have just proven in the previous lemma, that there is a way to construct $\mob_v$ for a drawing $z\in \Ce^Q$ for some graph $Q\in \obs^*$ by employing intersections of circles. This construction does not only give geometric meaning to the M\"obius mutation map, but will also enable us to prove Theorem \ref{th:miquelclifford} relating the Miquel and the Clifford construction, which in turn enables us to prove Theorem \ref{th:miqueldimers} on the connection between dimer statistics and Miquel dynamics.\\
In order to make progress later, we will study here a few more cross ratio identities of Clifford configurations. We borrow the following definition from Bobenko, Mercat and Suris \cite{nonlinearconformal}.
\begin{definition}
	Given a map $X: H_n \rightarrow \Ce$ from the hypercube $H_n$ to the plane, we say it is an \emph{integrable cross ratio system} if and only if for any $k, a,b\in\menge{1,2,\dots,n} , I\subset\menge{1,2,\dots,n}$:
	\begin{align}
		\shift{k}\cro(X_I, \shift{a} X_I, \shift{ab}X_I, \shift{b} X_I) = \cro(X_I, \shift{a} X_I, \shift{ab}X_I, \shift{b} X_I)
	\end{align}
	That is, on any 3-cube in $H_n$ the cross ratios on opposite 2-faces are the same.
\end{definition}
While this definition does not consider circle centers or intersection points, it can be applied to Clifford configurations and circle patterns.

\begin{lemma} \label{lem:cliffordintegrable}
	The $C_n$ configuration is an integrable cross ratio system.
\end{lemma}
\proof{
	A 3-cube in $C_n$ is simply $C_3$ (see figure \ref{fig:clifford}) and can be interpreted as an ideal hyperbolic tetrahedron in the Poincar\'e half-space model  (see \cite{hypgeom} for a reference on hyperbolic geometry). The circles represent the ideal boundary of the four involved planes and the intersection points are the four ideal vertices of the tetrahedron. The cross ratio $\cro(V,M_1,V_{12},M_2)$ has absolute value 1 and its argument is twice the intersection angles of the circles $c_1$ and $c_2$. Also, it is well known that intersection angles on opposite edges in ideal tetrahedra sum to $\pi$. Thus:
	\begin{align}
		&\cro(V,M_1,V_{12},M_2)\cro(V_{13},M_{123},V_{23},M_3) = 1\\
		\Rightarrow\qquad& \cro(V,M_1,V_{12},M_2) = \cro(M_3,V_{13},M_{123},V_{23}) = \shift{3} \cro(V,M_1,V_{12},M_2) 
	\end{align}
	Due to symmetry, this argument works for any cross ratios on pairs of opposing faces in $C_3$.\qed\\
}

\begin{lemma} \label{lem:tetrahedracrossratios}
	 In $C_4$ we have that
	\begin{align}
	\cro (V_I,\shift{ab} V_I,\shift{bc} V_I,\shift{ac} V_I) = \shift k \cro (V_I,\shift{ab} V_I,\shift{bc} V_I,\shift{ac} V_I)
	\end{align}
	for any $k\in \menge{1,2,3,4}$.	
\end{lemma}
\proof{We are comparing the cross ratios of the intersection points of two $C_3$ configurations situated in $C_4$. As already noted in Lemma \ref{lem:cliffordintegrable}, these $C_3$ configurations represent hyperbolic tetrahedra. In fact, these are two hyperbolic tetrahedra with the same intersection angles. It is well known \cite{hypgeom} that this is equivalent to the existence of a hyperbolic isometry, which is a M\"obius transformation of $\Ce$, that maps one of the tetrahedra onto the other. Thus, as M\"obius transformations do not change cross ratios, the lemma is proven.\qed\\}

Note that even though circle centers are not invariant under M\"obius transformations, the previous lemma does indeed relate cross ratios on intersection points of a $C_3$ configuration to the cross ratios of circle centers of that same $C_3$ configuration. This allows us to prove in the next lemma that star-ratios of circle centers are also preserved under the total shift $\Delta_{1234}$.

\begin{lemma} \label{lem:dualsr} Assume a $C_4$ configuration, where we picked the origin of the hypercube at the circle center $M$. The following relation holds:
\begin{align}
\sr(M; M_{12},M_{23},M_{34},M_{14}) = \shift{1234}\sr(M; M_{12},M_{23},M_{34},M_{14})
\end{align}
\end{lemma}
\proof{
	First, we check that this is true for the argument:
	\begin{align}
	\arg \sr(M; M_{12},M_{23},M_{34},M_{14}) &= \arg \cro (V_1,V_2,V_3,V_4)\\
	&= \arg \shift{1234} \cro (V_1,V_2,V_3,V_4)\\
	&= 	\arg \shift{1234} \sr(M; M_{12},M_{23},M_{34},M_{14})
	\end{align}	
	Thus it suffices to check that the quotient of the squares is equal to one. Call this quotient $S^2$ and decompose it into cross ratios of circle centers.
	\begin{align}
	S^2 :=& \left(\f{\sr(M; M_{12},M_{23},M_{34},M_{14})}{\sr(M_{1234}; M_{12},M_{23},M_{34},M_{14})} \right)^2
	\\=&\ \f{\cro(M,M_{12},M_{13},M_{23})}{\cro(M_{1234},M_{34},M_{13},M_{23})}\f{\cro(M,M_{34},M_{13},M_{14})}{\cro(M_{1234},M_{12},M_{13},M_{14})}
	\\&\cdot\f{\cro(M,M_{34},M_{24},M_{23})}{\cro(M_{1234},M_{34},M_{24},M_{14})}\f{\cro(M,M_{12},M_{24},M_{14})}{\cro(M_{1234},M_{12},M_{24},M_{23})}
	\end{align}			
	Notice that each cross ratio is a cross ratio of a $C_3$ configuration, therefore we are able to apply Lemma \ref{lem:tetrahedracrossratios}. We shift each cross ratio such that we arrive at an expression in which all the cross ratios cancel.
	\begin{align}
		S^2= \f{\cro(V_{123},V_{3},V_{2},V_{1})}{\cro(V_{124},V_{4},V_{1},V_{2})}\f{\cro(V_{134},V_{1},V_{4},V_{3})}{\cro(V_{234},V_{2},V_{3},V_{4})}
		\f{\cro(V_{234},V_{2},V_{3},V_{4})}{\cro(V_{123},V_{3},V_{2},V_{1})}\f{\cro(V_{124},V_{4},V_{1},V_{2})}{\cro(V_{134},V_{1},V_{4},V_{3})}
		=1
	\end{align}		
\qed \\}

For the special case that the circles at a quadrilateral are in a $C_4$ configuration, Lemma \ref{lem:dualsr} already proves the equivalence of Clifford and Miquel move. In particular, the so called integrable circle patterns \cite{nonlinearconformal} consist of $C_n$ configurations at all faces and vertices. As a result both the Miquel move and the Clifford move map integrable circle patterns to integrable circle patterns. However, we do not need this result here. Instead, we will prove in the next section \ref{sec:miquel} that locally, it is possible to reduce the case of general circle patterns to the case of an integrable cross ratio system.
 
% !TEX root = oupau.tex

\section{Circle patterns and Miquel dynamics}\label{sec:miquel}

In this section we prove the two main theorems by applying our previous results to the Miquel move. In Lemma \ref{lem:cphaverealsr} we have already proven that circle patterns feature real star-ratios, now we show that in the case of the sphere and the plane the reverse direction is true as well.

\begin{lemma} \label{lem:realsrleadtocp}
	Let $G\in \obs$ such that it is a decomposition of the sphere $S^2$ or the plane $\R^2$. Let $y\in \Ce^F$ be a face drawing of $G$ such that all star-ratios of $y$ are real. Then there exists a two parameter family of drawings $z\in \CP^G$ which are circle patterns such that $z^* = y$.
\end{lemma}
\proof{Let $v_0\in V_G$ be a fixed vertex, $z_0\in \Ce$ and set $z(v_0) = z_0$. We prove that we can extend this uniquely to a circle pattern on all of $G$. If $y=z^*$ is given and $z(v)$ is given, then $z(v')$ is uniquely determined for all $(v,v')\in E_G$. This is because $z(v')$ is the reflection of $z(v)$ about the line through the two incident circle centers. As $G$ is connected, there is at most one circle pattern such that $z(v_0) = z_0$. It remains to show that this extension always exists. For this it is sufficient to check that around any cycle in $G$, the propagation of the drawing via reflections closes. As all cycles in $G$ can be generated by face cycles, it suffices to check the closing condition around face cycles. This however is an immediate consequence of the identification of the total reflection angle with the argument of the star-ratio as in Lemma \ref{lem:cphaverealsr}.	
\qed\\}
Note that on surfaces like the torus, where cycles exist that are not generated by the face cycles, there need not exist a global circle pattern even if all the star-ratios are real. If it does exist, there will not necessarily be a one parameter family of circle patterns either.\\

\begin{figure}[th]
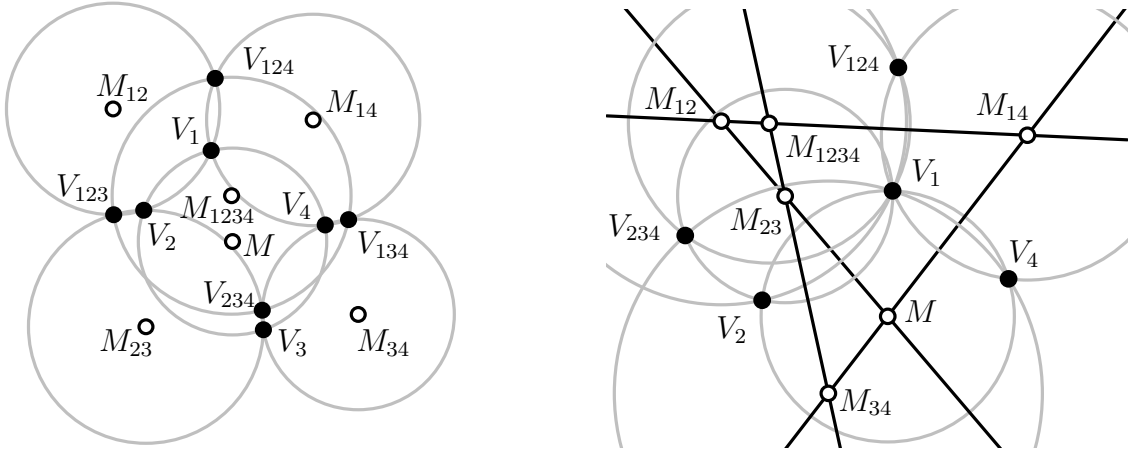
	
	\centering
	\begin{minipage}[b]{0.46\textwidth}
		\ggbmiquel
	\end{minipage}
	\begin{minipage}[b]{0.46\textwidth}
		\ggbmiquelmenelaus
	\end{minipage}
	
	\caption{The general Miquel configuration on the left and the Miquel-Menelaus configuration on the right, where triples of circle centers are on a line.}
	\label{fig:miquelmenelaus}
\end{figure}

In the next step, we show that the action of the Miquel move on the circle centers does not depend on a particular choice of circle pattern with these centers.
\begin{lemma} \label{lem:independentofradii} 
If we fix the six circle centers of a Miquel configuration (see figure \ref{fig:miquelmenelaus}), one of the circle intersection points can be chosen arbitrarily in $\Ce$. 
\end{lemma}
\proof{Choose the labeling of the Miquel configuration as in figure \ref{fig:miquelmenelaus}. We show that given $M,M_{12}, M_{23}, M_{34}, M_{14}$ then $M_{1234}$ is fixed no matter where we pick the intersection point $V_1$.
	Denote by $R_{AB}$ the reflection in the line through the points $A$ and $B$. Then all other intersection points arise as reflections of $V_1$ about lines through the given circle centers. In particular we have that:
	\begin{align}
		V_{123} = R_{M_{12}M_{23}}\circ R_{M_{12}M}\circ R_{M_{12}M_{14}}(V_{124})
	\end{align}
	Of course, this composition of three reflections about three lines intersecting in the point $M_{12}$ is a reflection about a line $\ell_{12}$ itself. This line $\ell_{12}$ is the perpendicular bisector of $V_{123}$ and $V_{124}$ which is therefore independent of the choice of $V_{1}$. By the same argument we find the lines $\ell_{23}, \ell_{34}$ and $\ell_{14}$ as well. Finally, $M_{1234}$ is the intersection point of these four lines, which is therefore also independent of the choice of $V_{1}$.\qed\\
	
}
Therefore the construction of the sixth circle center via Miquel dynamics does not require a realization by circles, knowing the five other circle centers is actually enough information. The next lemma is the last piece we need to prove our main theorems.

\begin{figure}[th]
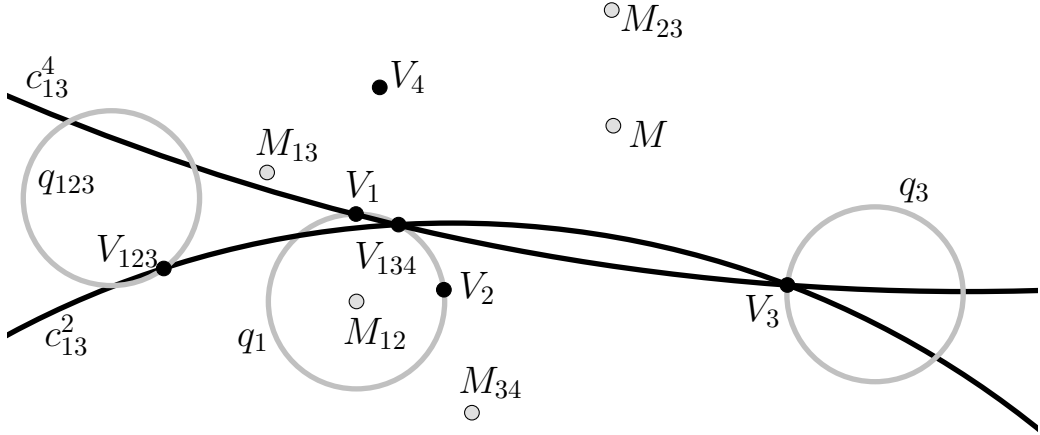
	
	\centering
	\ggblocalintegrability
	\caption{The Miquel configuration as we use it in the proof of Lemma \ref{lem:miquelsr}, case (ii).}
	\label{fig:localintegrability}
\end{figure}

\begin{lemma}
 \label{lem:miquelsr}
 Consider six circles in Miquel's configuration, then:
 \begin{align}
 	\sr(M; M_{12}, M_{23}, M_{34}, M_{14}) = \sr(M_{1234}; M_{12}, M_{23}, M_{34}, M_{14})
 \end{align}
\end{lemma}

\proof{We distinguish two cases:
\begin{itemize}
	\item[(i)] In the slightly degenerate case that $M$ is on the line $M_{k-1,k}M_{k,k+1}$, we will employ Menelaus's theorem again.
	\item[(ii)] If that is not the case, we show that there is a choice of circles such that the Miquel configuration is actually a part of a $C_4$ configuration.	
\end{itemize}
For the proof of case (i), let us relabel the indices such that $M\in M_{12}M_{23}$, see figure \ref{fig:miquelmenelaus}. The first consequence of $M\in M_{12}M_{23}$ is that $V_1=V_3$. The second consequence is that thus also $M\in M_{34}M_{14}$. We also note that by definition $V_{134} = V_{123} = V_1$. As a result, the perpendicular bisectors to $V_{124},V_{123}$ respectively $V_{124},V_{134}$ both coincide with the line $M_{12}M_{14}$. The same holds for the perpendicular bisectors to $V_{234},V_{123}$ respectively $V_{234},V_{134}$, which are both the line $M_{23}M_{34}$. Finally we have that $M_{1234} = M_{12}M_{14}\cap M_{23}M_{34}$ and we realize that the circle centers are in Menelaus's configuration, see figure \ref{fig:miquelmenelaus}.\\
In Menelaus's configuration we can use the same argument as in Lemma \ref{lem:cliffordmenelaus} to show that:
\begin{align}
\mr_1 &= \mr(M_{1234},M_{23},M_{34},M,M_{14},M_{12}) = -1\\
\mr_2 &= \mr(M_{1234},M_{34},M_{23},M,M_{12},M_{14}) = -1\\
\implies\qquad \f{\mr_1}{\mr_2} &= \f{\sr(M; M_{12},M_{23},M_{34},M_{14})}{\sr(M_{1234}; M_{12},M_{23},M_{34},M_{14})} = 1
\end{align}
For the proof of case (ii) we show that: Given the six midpoints $M,M_{12}, M_{23}, M_{34},M_{14}, M_{1234}$ of a Miquel configuration in general position, there is a choice of $V_1$ such that the configuration can be extended to a Clifford four circle configuration. \\
To show that such a choice exists, we introduce a parameter $t\in [0,1]$ moving $V_1(t)$ along a small circle $q_1$ centered at $M_{12}$ and then use a continuity argument, see also figure \ref{fig:localintegrability}. In particular, we will show that the circles $c_{13}^2(t)$ and $c_{13}^4(t)$ defined by the point triplets $V_1(t),V_3(t),V_{123}(t)$ and $V_1(t),V_3(t),V_{134}(t)$ coincide for some $t_0$ and therefore the choice $V_1(t_0)$ determines a $C_4$ configuration.\\
Because all the intersection points are determined by reflections about lines given by the circle centers, if $V_1(t)$ moves along a circle so do the other intersection points. We also define the circles $q_3,q_{123},q_{134}$ as the trajectories of the points $V_3(t),V_{123}(t),V_{134}(t)$. Because $V_{134}(t)$ is given by the successive reflections of $V_1(t)$ about the lines $M_{12}M_{23}$ and $M_{12}M$ the circle $q_{134}$ is in fact the same one as $q_1$. We choose the radius of the circle $c_{12}$, and therefore of all trajectory circles small enough such that the circles $q_1,q_3$ and $q_{123}$ are disjoint. It suffices to show that none of the centers $m_1,m_3,m_{123}$ of these circles coincide. Because $m_3$ is the reflection of $m_1=M_{12}$ about the line $MM_{34}$, these two could only coincide if $M,M_{12},M_{34}$ were collinear, which we have treated in case (i). Similarly, if $m_{123} = m_1$ or $m_{123} = m_3$, then $M,M_{14}, M_{34}$ would be collinear.\\
Now if we choose the circle radii very small, it is clear that of the two intersection points of $c_{13}^2(t)$ with $q_1$, one moves clockwise and one moves counter clockwise on $q_1$. Because $V_{134}(t)$ is a constant rotation away from $V_1(t)$ on $q_1$, both travel in the same direction on $q_1$. Therefore $V_{134}(t)$ will be on $c_{13}^2(t)$ for some $t$ which we can therefore select as $t_0$. Because $c_{13}^2(t_0)$ and $c_{13}^4(t_0)$ have the three points $V_1(t_0),V_3(t_0), V_{134}(t_0)$ in common they are in fact the same.\\
This proves that we can choose the intersection points such that $M,M_{12}, M_{23}, M_{34},M_{14}, M_{1234}$ are centers of a Clifford four circle configuration and thus it proves this case due to Lemma \ref{lem:dualsr}.
\qed\\
}

We have now shown that indeed the Miquel move on the circle centers acts as the M\"obius mutation does. As in the previous section, let us put this into a lemma.
\begin{lemma} \label{lem:miquelmobius}
		Let $G\in \obs$ be a graph, $z\in \Ce_0^G$ be a circle pattern and let $f\in F$ be a valid face with the four neighbours $f_1,f_2,f_3,f_4$, such that $z^*$ does not map the five points $f,f_1,f_2,f_3,f_3$ to a common line. The M\"obius move (see definition \ref{def:miquelmove}) coincides with the M\"obius mutation map on the induced face drawing $z^*$:
	\begin{align}
		(\miq_f(z))^* = \mob_f(z^*)
	\end{align}
\end{lemma}

As in the case of the Clifford move, it is possible to interpret the M\"obius mutation map as the extension of the Miquel move on the circle centers to the degenerate case that all the five centers are on a common line.\\
We can now proceed to prove the two main theorems. As a short reminder, Theorem \ref{th:miquelclifford} states that the Miquel and the Clifford move coincide and Theorem \ref{th:miqueldimers} states that the Miquel move behaves like urban renewal on the dimer model induced by the dual edge lengths.

\proof[P\textbf{roof of Theorem \ref{th:miquelclifford}}]{
	Let $z$ be a circle pattern of $G\in \obs$ and $f$ be a valid face.	In Lemma \ref{lem:cliffordmobius} and Lemma \ref{lem:miquelmobius} we have proven the following two equalities:
	\begin{align}
		\cli_f(z^*) &= \mob_f(z^*)\\
		(\miq_f(z))^* &= \mob_f(z^*)		
	\end{align}
	Together these two equations imply the theorem, namely that:
	\begin{align}
		\cli_f(z^*) = (\miq_f(z))^*
	\end{align}	
	%Let $G \in \obs$ be a graph and a $f\in F$ a degree four face. Also, let $z\in \CP^G$ be a circle pattern of $G$ and let $z^* \in \Ce^F$ be the dual drawing induced by the circle centers of $z$. By Lemma \ref{lem:miquelsr} the Miquel move $\miq_f$ and the Clifford move $\cli_{f^*}$ both leave the star-ratio invariant. However, because of the results of section \ref{sec:mobius} on the properties of the M\"obius mutation map, there is at most one such point different from $z^*(f)$. Therefore the Clifford and the Miquel move take $z^*(f)$ to the same point and the theorem is proven.\qed\\
	\qed\\
}

\proof[P\textbf{roof of Theorem \ref{th:miqueldimers}}]{
	In the case of a Kasteleyn circle pattern $z\in \CP^G$ all the star-ratios of $z^*\in \Ce^F$ are real and positive. Therefore in this case, the star-ratios only depend on the absolute values of the differences. Thus, if we consider the dimer edge weights $\psi$ on $G$ given by $z$ as in definition \ref{def:geometricdimer}, we conclude that:
	\begin{align}
		(\tau(\psi(z)))(v) = (\sr(z^*))(v^*) \qquad \forall v \in V
	\end{align}		
	Hence, the star-ratios are indeed the face weights of the associated dimer model. We have proven that the star-ratios transform under the Miquel move as they do under the M\"obius mutation map in Lemma \ref{lem:srmobius}. In the case that the star-ratios are real, these are indeed real formulae as well. We observe that they are exactly those of Lemma \ref{lem:urbanfaceweights} on the change of face weights $\tau$ in urban renewal. Thus the theorem is proven.\qed\\	
}

\section{Lattice dynamics} \label{sec:lattice}

Even though the new results of this paper are of local nature, we want to shortly outline the connection to dynamics on lattices. In particular, Miquel dynamics have originally been defined as an evolution of combinatorics on periodic $\Z^2$ lattices, and the study of Clifford lattices has taken place on the octahedral lattice, a sublattice of $\Z^3$.\\
Define the parity $|p|$ of a point $p = (p_1,p_2,p_3)\in \Z^3$ as follows:
\begin{align}
	|\cdot|: \Z^3 \rightarrow \Z_2,\qquad |p| = \sum_{k=1}^3p_k \mod 2
\end{align}

\begin{definition}
	The three dimensional \emph{octahedral lattice} $\Z_0^3$ consists of all points in $\Z^3$ that have even parity, together with edges between the nearest neighbours:
	\begin{align}
		V\left(\Z_0^3\right) &= \menge[p\in \Z^3]{|p| = 0}\\
		E\left(\Z_0^3\right) &= \menge[(p,p')\in V\left(\Z_0^3\right) \times V\left(\Z_0^3\right)]{\max_{k=1,2,3} |p_k-p_k'| = 1 }
	\end{align}
	Define a \emph{level slice} $S_k$ of the octahedral lattice as the vertices with third coordinate in $\menge{k,k+1}$ together with the edges of $\Z_0^3$ that connect two different layers.
	\begin{align}
		V(S_k) &= \menge[p\in V\left(\Z_0^3\right)]{p_3 \in \menge{k,k+1}}\\
		E(S_k) &= \menge[(p,p')\in E\left(\Z_0^3\right)]{p,p'\in V\left(\Z_0^3\right) \mbox{ and } p_3 \neq p'_3}
	\end{align}
\end{definition}

It is clear that such a level slice $S_k$ has the combinatorics of the $\Z^2$ lattice. As such, it can be bi-partitioned into the set of points which have $p_3 = k$ and the set of points which have $p_3=k+1$. 
\begin{align}
	S^0_k &:= \menge[p\in V\left(S_k\right)]{p_3 = k}\\
	S^+_k &:= \menge[p\in V\left(S_k\right)]{p_3 = k+1}
\end{align}
It follows that $S_k$ and $S_{k+1}$ are two $\Z^2$ lattices that as subsets of the octahedral lattice agree on every second point, that is $S^+_k= S^0_{k+1}$.\\
Set $G=S_k\in \obs^*$ considering it as a surface graph. Then the combinatoric 4-mutation at every vertex in $S^0_k$ gives a new graph $G'$ isomorphic to $\Z^2$ which we can identify with $S_{k+1}$, such that $S^+_k = S^0_{k+1}$ and all the vertices we mutated at are in $S^0_k$ before mutation and in $S^+_{k+1}$ after mutation. In that sense mutation replaces a vertex with the vertex two steps into the 3-direction of the octahedral lattice.\\
Let $z \in \Ce^{V_k}$ be a drawing where $V_k = V(S_k)$. Then by the combinatorial identification above and the M\"obius mutation map this drawing determines a unique function $z: V\left(\Z_0^3\right) \rightarrow \Ce$ from the octahedral lattice to the complex plane.
\begin{align}
	z(p_{x,y,z+1}) = \mob(z(p_{x+1,y+1,z}), z(p_{x+1,y-1,z}), z(p_{x-1,y-1,z}), z(p_{x-1,y+1,z}))(z(p_{x,y,z-1}))
\end{align}
Functions on octahedral lattices which are constructed as above have been characterized as discrete integrable systems in \cite{wolfgangclifford}. A system is called \emph{discretely integrable} if there exists a function on a lattice such that it satisfies a given set of local equations at each point in the lattice.\\
It is clear from the construction above that this is a 3-d integrable system in the sense that the data that determines the whole system (the Cauchy data) is 2-d, as we can choose freely the function on all points of a given slice in the octahedral lattice. In the case of circle patterns it is clear that we can choose as Cauchy data any circle pattern with the combinatorics of $\Z^2$. Such a circle pattern can be constructed row by row. Given a row, each circle center of the next row has to be on the bisector of the two corresponding intersection points of the given row, therefore there is a 1-dimensional real degree of freedom per circle. If we are interested in Kasteleyn circle patterns, then each circle center has to be chosen on an $\R^+$ ray and the Cauchy data is therefore still 2-d.\\
Given fixed star-ratios on $\Z^2 = S_k$, then Lemma \ref{lem:srmobius} asserts that we know the star-ratios on all slices $S_{k'}$ of the octahedral lattice. However, the knowledge of all the star-ratios on $\Z^2$ does not fix the drawing $z \in \Ce^{\Z_2}$. Instead, the fixed star-ratios determine a 2-d discrete integrable system themselves. Given the drawing $z$ on a strip $\Z \times \menge{l,l+1}$, the star-ratios determine the whole drawing on $\Z^2$. This strip is in fact the intersection of two level slices in two different directions. In this case, the prescription of real star-ratios corresponds to the case of circle patterns and the prescription of real and positive star-ratios corresponds to the case of Kasteleyn circle patterns.\\
An octahedral lattice constructed by the M\"obius mutation from given data on $\Z^2$ carries additional symmetry: The three space directions of the lattice are indistinguishable. This is a consequence of the symmetric nature of the M\"obius mutation map as explained in section \ref{sec:mobius}. In particular, if the star-ratios are the same on the two opposing stars in an octahedron of the lattice, then they are equal on all three stars in that octahedron. Moreover, if we prescribe real star-ratios in the given Cauchy data, then all the star-ratios of the whole lattice will be real, see Lemma \ref{lem:srmobius}. Therefore an octahedral lattice with real star-ratios induces three families of $\Z^2$ circle patterns, each family consisting of circle patterns related via Miquel dynamics. However, this symmetry does not hold for Kasteleyn circle patterns. If one level slice $S_k$ corresponds to a Kasteleyn circle pattern, then all the parallel slices are Kasteleyn as well, but the slices in the other directions are necessarily not Kasteleyn. This is an immediate consequence of the identities in equation \eqref{eq:sroctahedron}, which shows that positive star-ratios in one direction imply negative star-ratios in the other directions.
 
\section{Concluding remarks} \label{sec:conclusion}
Circle patterns are an inherently M\"obius invariant object, in the sense that circle patterns remain circle patterns after a M\"obius transformation. However, the object we used for the definition of the dimer model associated to a circle pattern are not. In particular M\"obius transformations do not map circle centers to circle centers, distances are not preserved and neither are star-ratios. It is therefore an interesting question how M\"obius transformations affect the statistics. It should be mentioned though that circle centers are closely related to Similarity geometry as a subgeometry of M\"obius geometry. If we choose one particular point  $P_\infty \in \Ce$, then we can define the so called conformal center of circle $c$ as reflection of $P_\infty$ about $c$. If $P_\infty = \infty$ then the conformal centers are just the usual circle centers. In this definition the star-ratio actually factors into cross ratios:
\begin{align}
	\sr(P_0; P_1,P_2,P_3,P_4) = \cro(P_1,P_0,P_2,P_\infty)\cro(P_3,P_0,P_4,P_\infty)
\end{align}
If we apply a M\"obius transformation to all the points including $P_\infty$, the star-ratios are invariants. However, this is really just invariance under Similarity transformations if we mod out transformations that do not fix $P_\infty$.\\

Another interesting question is if indeed Miquel dynamics are about circle patterns, or if they are slightly more fundamentally about real star-ratios and the existence of a consistent reflection system. By a consistent reflection system we mean the fact that the primal points are reflected about lines through the dual points. These two notions coincide on the sphere and the plane, but they differ on the torus. In fact, neither the equivalence to urban renewal nor to Clifford lattices depends on the existence of the global circle pattern, only on the existence of the reflection system.\\

Furthermore, we have shown that the dual star-ratios under Miquel dynamics carry the structure of a discrete integrable system as investigated in \cite{dimerintegrable} for graphs on the torus. However, it is unclear whether given star-ratios correspond to unique circle patterns or reflection systems on the torus. It is clear that on the infinite plane they are not sufficient, as discussed in section \ref{sec:lattice} where we showed that star-ratios only determine a $\Z^2$ lattice if one dimensional Cauchy-data are given.\\

Both the dimer and the Ising model have been studied on isoradial graphs and are known to be critical  \cite{isoradialdimercritical, baxter}. For the case of the Ising model there has been recent progress and it was proven by Lis \cite{criticalisingcirclepattern} that one can assign weights to a general circle pattern such that the resulting Ising model is critical. It seems a natural question to ask whether a similar result, using the edge weights for general circle patterns introduced here as a generalization of Kenyon's weights for isoradial graphs \cite{dimersisoradial}, give rise to a dimer model that is critical as well.\\

\section*{Funding}
This research was supported by the Deutsche Forschungsgemeinschaft (DFG) Collaborative Research Center TRR 109 ``Discretization in Geometry and Dynamics''.

\section*{Acknowledgments}
The author would like to thank Ananth Sridhar, Boris Springborn, Wolfgang Schief and Jan Techter for discussions and good advice.

\bibliographystyle{alpha}
\bibliography{references}

\end{document}